\documentclass[10pt]{amsart}
\usepackage{amscd}
\usepackage{amssymb}

\newtheorem{thm}[equation]{Theorem}
\newtheorem{cor}[equation]{Corollary}
\newtheorem{prop}[equation]{Proposition}
\newtheorem{lem}[equation]{Lemma}
\theoremstyle{definition}
\newtheorem{dfn}[equation]{Definition}
\newtheorem{rem}[equation]{Remark}
\newtheorem{exa}[equation]{Example}

\newtheorem{prob}[equation]{Problem}
\numberwithin{equation}{section}

\newcommand{\iso}{\stackrel{\simeq}{\rightarrow}}
\newcommand{\inj}{\hookrightarrow}
\newcommand{\surj}{\twoheadrightarrow}
\newcommand{\ar}{\rightarrow}
\newcommand{\opn}{\operatorname}
\newcommand{\cat}[1]{\operatorname{\mathsf{#1}}}

\newcommand{\blnk}[1]{\mbox{\hspace{#1}}}
\newcommand{\rmitem}[1]{\item[\text{\textup{(#1)}}]}

\newcommand{\mfrak}[1]{\mathfrak{#1}}
\newcommand{\mcal}[1]{\mathcal{#1}}
\newcommand{\msf}[1]{\mathsf{#1}}

\newcommand{\mrm}[1]{\mathrm{#1}}
\newcommand{\mbb}[1]{\mathbb{#1}}
\newcommand{\gfrac}[2]{\genfrac{[}{]}{0pt}{}{#1}{#2}}
\newcommand{\tup}[1]{\textup{#1}}

\title[Duality and Morita Equivalence]
{Dualizing Complexes, Morita Equivalence and the Derived Picard
Group of a Ring}
\author{Amnon Yekutieli}
\address{Department of Theoretical Mathematics,
The Weizmann Institute of Science,
Rehovot 76100, ISRAEL}
\date{30.8.98}
\email{amnon@wisdom.weizmann.ac.il}
\thanks{This work was partially supported by the US-Israel Binational
Science Foundation}

\begin{document}


\begin{abstract}
Two rings $A$ and $B$ are said to be derived Morita equivalent
if the derived categories
$\msf{D}^{\mrm{b}}(\cat{Mod} A)$ and
$\msf{D}^{\mrm{b}}(\cat{Mod} B)$
are equivalent. By results of Rickard in \cite{Ri1} and
\cite{Ri2}, if $A$ and $B$ are
derived Morita equivalent algebras over a field $k$, then
there is a complex of bimodules $T$ s.t.\ the functor
$T \otimes^{\mrm{L}}_{A} - : \msf{D}^{\mrm{b}}(\cat{Mod} A) \ar
\msf{D}^{\mrm{b}}(\cat{Mod} B)$
is an equivalence. The complex $T$ is called a tilting complex.

When $B = A$ the isomorphism classes of tilting complexes
$T$ form the derived Picard group $\opn{DPic}(A)$.
This group acts naturally on the Grothendieck group
$\opn{K}_{0}(A)$.

We prove that when the algebra $A$ is either local or commutative,
then any derived Morita equivalent algebra $B$ is actually Morita
equivalent. This enables us to compute $\opn{DPic}(A)$ in these
cases.

Assume $A$ is noetherian. Dualizing complexes over $A$ were
defined in \cite{Ye}. These are complexes of bimodules which
generalize the commutative definition of \cite{RD}.
We prove that the group $\opn{DPic}(A)$ classifies the set of
isomorphism classes of dualizing complexes.
We use this classification to deduce properties of rigid dualizing
complexes, as defined by Van den Bergh in \cite{VdB}.

Finally we consider finite $k$-algebras.
For the algebra $A$ of upper triangular $2 \times 2$ matrices over $k$,
we prove that $t^{3} = s$, where $t, s \in \opn{DPic}(A)$
are the classes of $A^{*} := \opn{Hom}_{k}(A, k)$ and
$A[1]$ respectively.
In the Appendix by Elena Kreines this result is generalized to
upper triangular $n \times n$ matrices, and it is shown that
the relation $t^{n + 1} = s^{n - 1}$ holds.
\end{abstract}

\maketitle

\setcounter{section}{-1}
\section{Introduction}

Let $A$ and $B$ be two rings. Recall that according to
Morita Theory, any equivalence between the categories of left modules
$\cat{Mod} A \ar \cat{Mod} B$ is realized by a $B$-$A$-bimodule
$P$, progenerator on both sides, as the functor
$M \mapsto P \otimes_{A} M$.

Happel \cite{Ha}, Cline, Parshall and Scott \cite{CPS} and
Rickard \cite{Ri1}, \cite{Ri2} generalized Morita theory to
derived categories.
Let $A$ and $B$ be algebras over a field $k$. Rickard proved that
if the derived categories
$\cat{D}^{\mrm{b}}(\cat{Mod} A)$ and
$\cat{D}^{\mrm{b}}(\cat{Mod} B)$
are equivalent, then there is a complex
$T \in \cat{D}^{\mrm{b}}(\cat{Mod} (B \otimes_{k} A^{\circ}))$
such that the functor
$T \otimes^{\mrm{L}}_{A} - :
\cat{D}^{\mrm{b}}(\cat{Mod} A) \ar
\cat{D}^{\mrm{b}}(\cat{Mod} B)$
is an equivalence. Here $A^{\circ}$ denotes the opposite algebra.
A complex $T$ with this property is called a {\em tilting complex},
and the algebras $A$ and $B$ are said to be {\em derived Morita
equivalent}.

In Section 1 we recall some facts on derived categories of bimodules
from \cite{Ye}. Then we reproduce Rickard's results in the
formulation needed for this paper. See Remark \ref{rem1.1}
regarding the generalization to an arbitrary commutative base
ring $k$.

In Section 2 we prove that if $A$ is either local or commutative
then any derived Morita equivalent algebra $B$ is actually
Morita equivalent (in the ordinary sense).
Specifically if
$T \in \cat{D}^{\mrm{b}}(\cat{Mod} (B \otimes_{k} A^{\circ}))$
is a tilting complex then $T \cong P[n]$ for some invertible bimodule
$P$ and some integer $n$ (in the commutative case $\opn{Spec} A$
is assumed connected). See Theorems \ref{thm2.7} and \ref{thm2.4}.

When $A = B$ the isomorphism classes of tilting complexes
form a group,
called the {\em derived Picard group} $\mrm{DPic}(A)$.
The operation is
$(T_{1}, T_{2}) \mapsto T_{1} \otimes^{\mrm{L}}_{A} T_{2}$,
the identity is $A$ and the inverse is
$T \mapsto T^{\vee} := \mrm{R} \opn{Hom}_{A}(T, A)$.
Let $s \in \mrm{DPic}(A)$ be the class of the complex
$A[1]$. Then the subgroup $\langle s \rangle$ is isomorphic to
$\mbb{Z}$.
When $A$ is local we show that
$\mrm{DPic}(A) \cong \mbb{Z} \times \opn{Out}_{k}(A)$,
where $\opn{Out}_{k}(A)$ denotes the group of outer $k$-algebra
automorphisms (see Proposition \ref{prop3.4}).
When $A$ is commutative then
$\mrm{DPic}(A) \cong \mbb{Z}^{m} \times
\opn{Aut}_{k}(A) \ltimes \opn{Pic}_{A}(A)$,
where $m$ is the number of connected components of $\opn{Spec} A$
and $\opn{Pic}_{A}(A)$ is the usual commutative Picard group
(Proposition \ref {prop3.3}).
If $A$ is noetherian let
$\opn{K}_{0}(A) = \opn{K}_{0}(\cat{Mod}_{\mrm{f}}(A))$
be the Grothendieck group. Then there is a representation
$\chi_{0} : \mrm{DPic}(A) \ar \opn{Aut}(\opn{K}_{0}(A))$,
with $\chi_{0}(s) = -1$.

In Section 4 we suppose $A$ is noetherian. Then we have
the notion of {\em dualizing complex}
$R \in \cat{D}^{\mrm{b}}(\cat{Mod} A^{\mrm{e}})$,
where $A^{\mrm{e}} := A \otimes_{k} A^{\circ}$
(Definition \ref{dfn4.1}).
Dualizing complexes over noncommutative algebras were introduced
in \cite{Ye}, generalizing the commutative definition of \cite{RD}.
Unlike the commutative case, where any two dualizing complexes
$R_{1}, R_{2}$ satisfy
$R_{2} \cong L[n] \otimes_{A} R_{1}$ with $L$ an invertible module
and $n$ an integer, when $A$ is noncommutative
there is no such uniqueness.
The question arose how to classify all isomorphism classes.
We prove in Theorem \ref{thm4.1}
that given a dualizing complex $R_{1}$, any other complex
$R_{2}$ is dualizing if{f} 
$R_{2} \cong T \otimes^{\mrm{L}}_{A} R_{1}$
for some tilting complex $T$. Moreover this $T$ is
unique up to isomorphism. Therefore the group $\opn{DPic}(A)$
classifies the isomorphism classes of dualizing complexes.

Next, in Section 5,
we consider {\em rigid dualizing complexes}, which were defined by
Van den Bergh \cite{VdB}.
One of his results was that a rigid dualizing complex is unique up to
an isomorphism in $\cat{D}(\cat{Mod} A^{\mrm{e}})$.
We prove that this isomorphism is unique
(Theorem \ref{thm5.1}).
If $A$ is finitely generated as $k$-algebra and finite over its center
then it has a rigid dualizing complex
(Proposition \ref{prop5.1}).
If $A$ is Gorenstein and has a rigid dualizing $R$ complex then
$R$ is also a tilting complex, and
$R^{\vee} \cong \mrm{R} \opn{Hom}_{A^{\mrm{e}}}(A, A^{\mrm{e}})$
(Proposition \ref{prop5.5}). This also generalizes a result of
Van den Bergh.

In Section 6 we look at a finite $k$-algebra $A$.
The bimodule $A^{*} := \opn{Hom}_{k}(A, k)$ is the rigid dualizing
complex of $A$. If $A$ is a Gorenstein algebra then $A^{*}$ is also a
tilting complex, in which case we denote its class in $\opn{DPic}(A)$
by $t$. If moreover $A$ has finite global dimension then
$\chi_{0}(t) = -c$, where $c$ is the Coxeter transformation of
\cite{ARS} Chapter VIII.

Finally in Proposition \ref{prop6.3}
we examine the group $\opn{DPic}(A)$ for the algebra
$A =
\left[ \begin{smallmatrix}
k & k \\ 0 & k \end{smallmatrix} \right]$.
Note that this is the smallest $k$-algebra which is neither
commutative nor local. The ordinary noncommutative Picard
group $\mrm{Pic}(A)$ is trivial here. On the other hand, we
prove that $t^{3} = s$, so
$\mrm{DPic}(A) \neq \langle s \rangle$.
In the Appendix by Elena Kreines the calculation is carried out for
an $n \times n$ upper triangular matrix algebra, $n \geq 2$. She
proves that $t^{n + 1} = s^{n - 1}$ in this case.

\medskip \noindent \textbf{Acknowledgments.}\
Work on the paper started in conversations with K.\ Ajitabh, during
a visit funded by a grant from the US-Israel Binational Science
Foundation. I am grateful to M.\ Van den Bergh, V.\ Hinich,
I.\ Reiten, J.\ Rickard and B. Keller for their helpful comments; 
to E. Kreines who contributed the Appendix; and to the referee 
for finding a delicate mistake in the original formulation of
Lemma \ref{lem2.5}.

\section{Morita  Equivalence}

Let $k$ be a fixed base field.
All $k$-algebras will be associative with $1$.
Given a $k$-algebra $A$ we denote by $A^{\circ}$
the opposite algebra, and by $A^{\mrm{e}}$ the enveloping algebra
$A \otimes A^{\circ}$ (where $\otimes = \otimes_{k}$ throughout).
Our modules will be by default left modules,
and with this convention an $A^{\circ}$-module
will mean a right $A$-module.
Given another $k$-algebra $B$, an $(A \otimes B^{\circ})$-module $M$ is
then just an $A$-$B$-bimodule ${}_{A}M_{B}$, central over $k$.

We write $\cat{Mod} A$ for the category of $A$-modules.
Let $\cat{D}(\cat{Mod} A)$ be the derived category of $A$-modules,
and for $\star = -, +, \mrm{b}$ let
$\cat{D}^{\star}(\cat{Mod} A)$ be the appropriate full subcategory
(conventions as in \cite{RD}).

The forgetful functor
$\cat{Mod}(A \otimes B^{\circ}) \ar \cat{Mod} A$
is exact and so induces a functor
$\cat{D}^{\star}(\cat{Mod}(A \otimes B^{\circ})) \ar
\cat{D}^{\star}(\cat{Mod} A)$.
Now $A \otimes B^{\circ}$ is a projective $A$-module, so
any projective (resp.\ flat, injective) $(A \otimes B^{\circ})$-module
is projective (resp.\ flat, injective) over $A$.

Consider $k$-algebras $A, B, C$. For complexes
$M \in \cat{D}(\cat{Mod}(A \otimes B^{\circ}))$ and
$N \in \cat{D}(\cat{Mod}(A \otimes C^{\circ}))$,
with either
$M \in \msf{D}^{-}$ or $N \in \msf{D}^{+}$,
there is a derived functor
\[ \mrm{R} \opn{Hom}_{A}(M, N) \in
\cat{D}(\cat{Mod}(B \otimes C^{\circ})) . \]
It is calculated by replacing $M$ with an isomorphic complex
in $\cat{D}^{-}(\cat{Mod}(A \otimes B^{\circ}))$
which consists of projective modules over $A$; or
by replacing $N$ with an isomorphic complex
in $\cat{D}^{+}(\cat{Mod}(A \otimes C^{\circ}))$
which consists of injective modules over $A$.

Likewise, for  complexes
$M \in \cat{D}^{-}(\cat{Mod}(B \otimes A^{\circ}))$
and
$N \in \cat{D}^{-}(\cat{Mod}(A \otimes C^{\circ}))$
there is a derived functor
\[ M \otimes^{\mrm{L}}_{A} N \in
\cat{D}(\cat{Mod}(B \otimes C^{\circ})) . \]
It is calculated by replacing $M$ with an isomorphic complex
in $\cat{D}^{-}(\cat{Mod}(B \otimes A^{\circ}))$
which consists of flat modules over $A^{\circ}$, or by doing the
corresponding thing for $N$.
In case $M$ has finite Tor dimension over $A^{\circ}$, i.e.\ it is
isomorphic in $\cat{D}(\cat{Mod} A^{\circ})$
to a bounded complex of flat $A^{\circ}$-modules, then
$M \otimes^{\mrm{L}}_{A} N$ is defined for an unbounded $N$ (and vice
versa).
For full details see \cite{Ye}.

Because the forgetful functors
$\cat{Mod}(A \otimes B^{\circ}) \ar \cat{Mod} A$
etc.\ commute with \linebreak
$\mrm{R} \opn{Hom}_{A}(-,-)$ and $- \otimes^{\mrm{L}}_{A} -$
there is no need to mention them explicitly.

Recall that a complex $M \in \cat{D}(\cat{Mod} A)$
is called {\em perfect} if it is isomorphic to a bounded complex of
finitely generated projective modules.
The full subcategory
$\cat{D}(\cat{Mod} A)_{\mrm{perf}} \subset
\cat{D}(\cat{Mod} A)$ consisting of perfect complexes is
triangulated, and the identity functors
\[ \cat{K}^{\mrm{b}}(\cat{proj} A) \ar
\cat{D}^{\mrm{b}}(\cat{Mod} A)_{\mrm{perf}} \ar
\cat{D}(\cat{Mod} A)_{\mrm{perf}} \]
are equivalences, where $\cat{proj} A$ is the additive category
of finitely generated projective $A$-modules.

\begin{lem} \label{lem1.2}
\begin{enumerate}
\item Suppose
$M \cong M_{1} \oplus M_{2} \in \cat{D}(\cat{Mod} A)$.
Then $M$ is perfect if{f} both $M_{1}$ and $M_{2}$ are perfect.
\item Let $M \in \cat{D}(\cat{Mod} A)$ be a perfect complex and $n$
an integer. If $\mrm{H}^{p} M = 0$ for all $p > n$ then
$\mrm{H}^{n} M$ is a finitely generated module.
\end{enumerate}
\end{lem}

\begin{proof}
1.\ See \cite{SGA6} Expos\'{e} I Proposition 4.17. 
\newline
2.\ Let $P \cong M$ where $P$ is a bounded complex of finitely 
generated projectives. Say $P^{p} = 0$ for $p > m \geq n$. 
By splitting $P^{n} \ar \cdots \ar P^{m}$ we obtain a surjection
$P^{n} \surj \mrm{H}^{n} M$.
\end{proof}

Given a complex $M \in \cat{D}(\cat{Mod} A)$, denote by
$\cat{add} M \subset \cat{D}(\cat{Mod} A)$
the class of all direct summands of finite direct sums of $M$.

We say a class $\cat{D}_{0}$ of objects of a triangulated category
$\cat{D}$ generates it if there is no triangulated subcategory
$\cat{D}'$, closed under isomorphisms, with
$\cat{D}_{0} \subset \cat{D}' \subsetneqq \cat{D}$.

At this point we wish to remind the reader of the classical Morita
Theory.

\begin{thm} \tup{(Morita Theory)}\ \label{thm1.2}
Let $A$ and $B$ be rings. Then the following are equivalent\tup{:}
\begin{enumerate}
\rmitem{i} The abelian categories $\cat{Mod} A$ and $\cat{Mod} B$ are
equivalent.
\rmitem{ii} There is a $B$-$A$-bimodule $P$, progenerator over $B$,
such that the canonical ring homomorphism
$A^{\circ} \ar \opn{End}_{B}(P)$
is bijective.
\rmitem{iii} There is a $B$-$A$-bimodule $P$ and an
$A$-$B$-bimodule $Q$ such that
$P \otimes_{A} Q \cong B$ and $Q \otimes_{B} P \cong A$ as
bimodules.
\end{enumerate}
\end{thm}

If $F : \cat{Mod} A \ar \cat{Mod} B$ is the equivalence, then
$P = F A$, $F M = P \otimes_{A} M$ and
$Q = \opn{Hom}_{B}(P, B)$.
In this case we say that $A$ and $B$ are Morita equivalent, and
we call a bimodule $P$ as above an {\em invertible $B$-$A$-bimodule}.

Following Rickard we make the following definition:

\begin{dfn}
Let $A$ and $B$ be rings.
If there is an equivalence of triangulated categories
$F : \cat{D}^{\mrm{b}}(\cat{Mod} A) \ar
\cat{D}^{\mrm{b}}(\cat{Mod} B)$
we say that $A$ and $B$ are {\em derived Morita equivalent}.
\end{dfn}

The generalization to complexes of the notion of invertible bimodule
is:

\begin{dfn} \label{dfn1.1}
Let $A, B$ be $k$-algebras, and let
$T \in  \cat{D}^{\mrm{b}}(\cat{Mod}(B \otimes A^{\circ}))$.
Suppose that:
\begin{enumerate}
\rmitem{i} $T \in \cat{D}^{\mrm{b}}(\cat{Mod} B)$
is a perfect complex, and $\cat{add} T$ generates
$\cat{D}^{\mrm{b}}(\cat{Mod} B)_{\mrm{perf}}$.
\rmitem{ii} The canonical morphism
$A \ar  \mrm{R} \opn{Hom}_{B}(T, T)$ in
$\cat{D}(\cat{Mod} A^{\mrm{e}})$
is an isomorphism.
\end{enumerate}
Then we call $T$ a {\em tilting complex}.
\end{dfn}

In \cite{Ri2} the name two-sided tilting complex was used.

\begin{exa}
In the notation of Theorem \ref{thm1.2}, if $P$ is $k$-central then
$P \in \cat{Mod} (B \otimes A^{\circ})$
is a tilting complex.
\end{exa}

The next theorem is an immediate consequence of \cite{Ri1}
Theorem 6.4, and \cite{Ri2} Theorem 3.3 and Proposition 4.1.
For the convenience of the reader we have included our own proof.

\begin{thm} \label{thm1.1} \tup{(Rickard)}\
The following conditions are equivalent for a complex
$T \in \cat{D}^{\mrm{b}}(\cat{Mod} (B \otimes A^{\circ}))$\tup{:}
\begin{enumerate}
\rmitem{i} The functor
\[ T \otimes^{\mrm{L}}_{A} - :
\cat{D}^{-}(\cat{Mod} A) \ar \cat{D}^{-}(\cat{Mod} B) \]
is an equivalence of triangulated categories.
\rmitem{i$'$} The functor $ T \otimes^{\mrm{L}}_{A} -$
preserves bounded complexes, and induces an equivalence  
of triangulated categories
\[ T \otimes^{\mrm{L}}_{A} - :
\cat{D}^{\mrm{b}}(\cat{Mod} A) \ar 
\cat{D}^{\mrm{b}}(\cat{Mod} B) . \]
\rmitem{ii} $T$ is a tilting complex.
\rmitem{iii} There exist a complex
$T^{\vee} \in \cat{D}^{\mrm{b}}(\cat{Mod}(A \otimes B^{\circ}))$
and isomorphisms
\[ \begin{aligned}
T^{\vee} \otimes_{B}^{\mrm{L}} T & \cong A \in
\cat{D}(\cat{Mod} A^{\mrm{e}}) , \\[1mm]
T \otimes_{A}^{\mrm{L}} T^{\vee} & \cong B \in
\cat{D}(\cat{Mod} B^{\mrm{e}}) .
\end{aligned} \]
\end{enumerate}
\end{thm}

\begin{proof}
(i) $\Rightarrow$ (i$'$):
First note that the identity functor
$\msf{K}^{-}(\cat{Proj} B) \ar \msf{D}^{-}(\cat{Mod} B)$
is an equivalence, where $\cat{Proj} B$
is the additive category of projective $B$-modules.
Now use \cite{Ri1} Proposition 6.1.

\medskip \noindent (i$'$) $\Rightarrow$ (ii):
Let
$F := T \otimes_{A}^{\mrm{L}} -$.
By \cite{Ri1} Propositions 6.2 and 6.3 we see that
$F$ restricts to an equivalence
$\cat{D}^{\mrm{b}}(\cat{Mod} A)_{\mrm{perf}} \ar
\cat{D}^{\mrm{b}}(\cat{Mod} B)_{\mrm{perf}}$.
Since $\cat{add} A$ generates
$\cat{D}^{\mrm{b}}(\cat{Mod} A)_{\mrm{perf}}$
it follows that $T = F A$ generates
$\cat{D}^{\mrm{b}}(\cat{Mod} B)_{\mrm{perf}}$.
Also $F$ induces isomorphisms
\[ \opn{Hom}_{ \cat{D}(\cat{Mod} A) }(A, A[i]) \iso
\opn{Hom}_{ \cat{D}(\cat{Mod} B) }(T, T[i]) , \]
so condition (ii) of the Definition \ref{dfn1.1} holds.

\medskip \noindent (ii) $\Rightarrow$ (i):
Let
\[ T^{\vee} := \mrm{R} \opn{Hom}_{B}(T, B) \in 
\cat{D}(\cat{Mod}(A \otimes B^{\circ})) . \]
Since
$T \in \cat{D}^{\mrm{b}}(\cat{Mod} B)_{\mrm{perf}}$
it follows that $T^{\vee}$ is bounded. 
Let
\[ F^{\vee} := T^{\vee} \otimes_{B}^{\mrm{L}} - :
\msf{D}^{-}(\cat{Mod} B) \ar \msf{D}^{-}(\cat{Mod} A) , \]
so
$F^{\vee} \cong \mrm{R} \opn{Hom}_{B}(T, -)$.
Now for any
$L \in \msf{D}^{-}(\cat{Mod} (B \otimes A^{\circ}))$,
$M \in \msf{D}^{-}(\cat{Mod} A)$
and
$N \in \msf{D}(\cat{Mod} B)$
one has an isomorphism
\[ \mrm{R} \opn{Hom}_{B}(L \otimes_{A}^{\mrm{L}} M, N) \cong
\mrm{R} \opn{Hom}_{A}(M, \mrm{R} \opn{Hom}_{B}(L, N)), \]
as can be seen by taking $L$ to be a complex of
$(B \otimes A^{\circ})$-projectives and
$M$ a complex of $A$-projectives. Therefore $F^{\vee}$
is a right adjoint to $F$, and condition (ii) of Definition 
\ref{dfn1.1} says that
$1_{\msf{D}^{-}(\cat{Mod} A)} \cong F^{\vee} F$.
Given any $M \in \msf{D}^{-}(\cat{Mod} B)$,
let $N$ be the cone on $F F^{\vee} M \ar M$.
Because
$F^{\vee} F F^{\vee} M \cong F^{\vee} M$
we find that
$\mrm{R} \opn{Hom}_{B}(T, N) \cong F^{\vee} N = 0$.
Now $\cat{add} T$ generates
$\cat{D}^{\mrm{b}}(\cat{Mod} B)_{\mrm{perf}}$
and $B \in \cat{D}^{\mrm{b}}(\cat{Mod} B)_{\mrm{perf}}$.
This implies $N = 0$ and hence
$F F^{\vee} \cong 1_{\msf{D}^{-}(\cat{Mod} B)}$
(cf.\ \cite{Ri1} Proposition 5.4).

\medskip \noindent (iii) $\Rightarrow$ (i):
The associativity of $- \otimes^{\mrm{L}} -$ implies that
$F F^{\vee} \cong 1_{\msf{D}^{-}(\cat{Mod} B)}$
and
$F^{\vee} F \cong 1_{\msf{D}^{-}(\cat{Mod} A)}$.

\medskip \noindent (ii) $\Rightarrow$ (iii):
Since $T$ is a tilting complex we have
\[ A \cong \mrm{R} \opn{Hom}_{B}(T, T) \cong
T^{\vee} \otimes_{B}^{\mrm{L}} T 
\in \msf{D}(\cat{Mod} A^{\mrm{e}}) . \]
By the proof of ``(ii) $\Rightarrow$ (i)'' the functor 
$F^{\vee} = T^{\vee} \otimes_{B}^{\mrm{L}} -$
is an equivalence, hence by ``(i) $\Rightarrow$ (ii)'',
$T^{\vee}$ is a tilting complex.
Writing
$T^{\vee \vee} := \mrm{R} \opn{Hom}_{A}(T^{\vee}, A)$,
the previous arguments show that 
$T^{\vee \vee} \otimes_{A}^{\mrm{L}} T^{\vee} \cong B
\in \msf{D}(\cat{Mod} B^{\mrm{e}})$.
But 
$T^{\vee \vee} \cong T^{\vee \vee} \otimes_{A}^{\mrm{L}}
T^{\vee} \otimes_{B}^{\mrm{L}} T \cong T$, 
and this completes the circle of the proof.
\end{proof}

\begin{cor} \label{cor1.1}
Let $A, B, C$ be $k$-algebras and let
$T \in \cat{D}^{\mrm{b}}(\cat{Mod} (B \otimes A^{\circ}))$
and
$S \in \cat{D}^{\mrm{b}}(\cat{Mod} (C \otimes B^{\circ}))$
be tilting complexes. 
\begin{enumerate}
\item $T \in \cat{D}^{\mrm{b}}(\cat{Mod} (A^{\circ} \otimes B))$
is a tilting complex, i.e.\ the roles of the algebras $A$ and $B$ 
in Definition \tup{\ref{dfn1.1}} can be exchanged. 
\item 
$T^{\vee} \in \cat{D}^{\mrm{b}}(\cat{Mod} (A \otimes B^{\circ}))$
from condition \tup{(iii)} of the theorem is a tilting complex, 
and it is unique up to isomorphism.
\item $S \otimes^{\mrm{L}}_{B} T \in
\cat{D}^{\mrm{b}}(\cat{Mod} (C \otimes A^{\circ}))$
is a tilting complex.
\item There are equivalences of triangulated categories
\[ \begin{aligned}
T \otimes_{A}^{\mrm{L}} - :\
& \cat{D}^{\star}(\cat{Mod} (A \otimes C^{\circ})) \ar
\cat{D}^{\star}(\cat{Mod} (B \otimes C^{\circ})) , \\
- \otimes_{B}^{\mrm{L}} T :\
& \cat{D}^{\star}(\cat{Mod} (C \otimes B^{\circ})) \ar
\cat{D}^{\star}(\cat{Mod} (C \otimes A^{\circ}))
\end{aligned} \]
with $\star = \mrm{b}, +, -, \emptyset$.
\end{enumerate}
\end{cor}

\begin{proof}
1, 2, 3. These are immediate consequences of condition (iii) of the 
theorem. \newline
4. Since $T$ is perfect over $A^{\circ}$ and over $B$ the functors 
are defined on the unbounded categories $\cat{D}$, and preserve
$\cat{D}^{\star}$. By way-out reasons (cf.\ \cite{RD} Proposition 
I.7.1(iv)) they are equivalences.
\end{proof}

We call the complex $T^{\vee}$ above {\em the inverse of} $T$.

The next theorem was shown to the author by V.\ Hinich.

\begin{thm} \label{thm1.5}
Let $B$ be a $k$-algebra, let
$T' \in \cat{D}^{-}(\cat{Mod} B)$ be a complex and let
$A := \opn{End}_{\cat{D}(\cat{Mod} B)}(T')$.
Assume
$\opn{Hom}_{\cat{D}(\cat{Mod} B)}(T', T'[i]) = 0$
for $i < 0$.
Then there is a complex
$T \in \cat{D}^{-}(\cat{Mod} (B \otimes A))$
s.t.\
$T \cong T'$ in $\cat{D}(\cat{Mod} B)$,
and the ring homomorphism
$A \ar \opn{End}_{\cat{D}(\cat{Mod} B)}(T)$
induced by the $A$-module structure of $T$ is bijective.
\end{thm}

\begin{proof}
We shall use ideas from homotopical algebra.
Suppose $C$ is a DGA (differential graded  algebra) over $k$, and
denote by $\cat{DGMod} C$ the category of DG $C$-modules. According to
\cite{Hi} Section 3, $\cat{DGMod} C$ is a closed model category in
the sense of Quillen \cite{Qu}. The weak equivalences in
$\cat{DGMod} C$ are the quasi-isomorphisms. Let
$\cat{D}(\cat{DGMod} C) = \cat{Ho}(\cat{DGMod} C)$
be the homotopy category, obtained by
inverting the weak equivalences. It is a triangulated category.
If $C$ is just a $k$-algebra (i.e.\ $C^{i} = 0$ for $i \neq 0$)
then
$\cat{DGMod} C = \cat{C}(\cat{Mod} C)$
and
$\cat{D}(\cat{DGMod} C) = \cat{D}(\cat{Mod} C)$.

According to \cite{Hi} Theorem 3.3.1
(or \cite{Ke2} Theorem 8.2), if
$C' \ar C$ is a quasi-isomorphism of DGAs, then the functor
$\cat{D}(\cat{DGMod} C) \ar \cat{D}(\cat{DGMod} C')$
gotten by restriction of scalars is an equivalence.

Given our complex $T'$, we may assume it consists of projective
$B$-modules. Define
$A'' := \opn{End}_{B}(T')$, which is a DGA,
and $A = \mrm{H}^{0} A''$. Let $A'$ be the truncation
$\sigma_{\leq 0} A''$, that is
\[ A' := \left(\cdots \ar {A''}^{-1} \ar
\opn{Ker}({A''}^{0} \ar {A''}^{1}) \ar 0 \ar \cdots \right) . \]
Since $A' \ar A''$ is a DGA homomorphism, we have
$T' \in \cat{DGMod} (B \otimes A')$.
On the other hand $A' \ar A$ is a quasi-isomorphism, and hence so
is
$B \otimes A'\ar B \otimes A$.
Consider the commutative diagram
\[ \begin{CD}
\cat{D}(\cat{Mod} (B \otimes A)) @>G >>
\cat{D}(\cat{DGMod} (B \otimes A')) \\
@VVV @VVV \\
\cat{D}(\cat{Mod} B) @> = >>
\cat{D}(\cat{Mod} B)
\end{CD} \]
where all the arrows are restriction of scalars.
Since $G$ is an equivalence, we can find a complex
$T \in \cat{D}(\cat{Mod} (B \otimes A))$
s.t.\
$G T \cong T'$ in $\cat{D}(\cat{DGMod} (B \otimes A'))$.
We may assume (by truncation) that
$T \in \cat{D}^{-}(\cat{Mod} (B \otimes A))$,
and then it has the desired properties.
\end{proof}

The following corollary is \cite{Ri2} Corollary 3.5. Our proof
is almost identical to B. Keller's in \cite{Ke1}.

\begin{cor} \label{cor1.6} \tup{(Rickard)}\
Let $A$ and $B$ be $k$-algebras, and let
$F : \cat{D}^{\mrm{b}}(\cat{Mod} A) \ar 
\cat{D}^{\mrm{b}}(\cat{Mod} B)$
be an equivalence of triangulated categories.
Then there exists a tilting complex
$T \in \cat{D}^{\mrm{b}}(\cat{Mod} (B \otimes A^{\circ}))$
with
$T \cong F A$ in $\cat{D}(\cat{Mod} B)$.
\end{cor}

\begin{proof}
According to \cite{Ri1} Propositions 6.1-6.3, $F$ restricts 
to an equivalence
$F : \cat{D}^{\mrm{b}}(\cat{Mod} A)_{\mrm{perf}}
\ar \cat{D}^{\mrm{b}}(\cat{Mod} B)_{\mrm{perf}}$.
Then $T' := F A \in \cat{D}(\cat{Mod} B)$
is a perfect complex,
$\cat{add} T'$ generates $\cat{D}^{\mrm{b}}(\cat{Mod} B)_{\mrm{perf}}$,
$\opn{End}_{\cat{D}(\cat{Mod} B)}(T') \cong A^{\circ}$
and \newline
$\opn{Hom}_{\cat{D}(\cat{Mod} B)}(T', T'[i]) = 0$
for $i \neq 0$. Now use Theorem \ref{thm1.5}.
\end{proof}

\begin{rem}
We did not assume that $F$ is $k$-linear in the corollary. But 
even when $F$ is $k$-linear, it is not known whether
the two functors $F$ and $T \otimes^{\mrm{L}}_{A} -$
are necessarily isomorphic. Rickard calls an equivalence of the form 
$T \otimes^{\mrm{L}}_{A} -$ {\em standard} (see \cite{Ri1} Section 
7 and \cite{Ri2} Definition 3.4).
\end{rem}

To finish off this section, consider a noetherian algebra $A$.
Then $\cat{Mod}_{\mrm{f}} A$, the category of
finitely generated modules, is abelian, and the category
$\msf{D}_{\mrm{f}}(\cat{Mod} A)$ of complexes with
finitely generated cohomologies is triangulated.

\begin{prop} \label{prop1.6}
If $A$ and $B$ are both noetherian and
$T \in \msf{D}^{\mrm{b}}(\cat{Mod} (B \otimes A^{\circ}))$
is a tilting complex then
\[ T \otimes^{\mrm{L}}_{A} - :
\msf{D}^{\star}_{\mrm{f}}(\cat{Mod} A) \ar
\msf{D}^{\star}_{\mrm{f}}(\cat{Mod} B) \]
is an equivalence of triangulated categories for
$\star = \mrm{b}, +, -, \emptyset$.
\end{prop}

\begin{proof}
Since $T \otimes^{\mrm{L}}_{A} -$ is a way-out functor in both
directions and
$T \otimes^{\mrm{L}}_{A} A = T \in
\msf{D}^{\star}_{\mrm{f}}(\cat{Mod} B)$
the proposition follows from \cite{RD} Proposition I.7.3.
\end{proof}

\begin{rem} \label{rem1.1}
Throughout the paper the base ring $k$ is a field. But it is easy to 
see that everything in Sections 1-3 will remain valid if we let 
$k$ be an arbitrary commutative ring, as long as the $k$-algebras 
$A, B, C$ are assumed to be {\em projective} $k$-modules. With a 
mild modification of the proofs we can even assume these algebras 
are only {\em flat} $k$-modules.
 
For the general situation here is an approach suggested by V. Hinich 
and B. Keller. Consider a DG $k$-algebra 
$B \otimes_{k}^{\mrm{L}} A^{\circ} =
\tilde{B} \otimes_{k} \tilde{A}^{\circ}$,
where $\tilde{A}$, $\tilde{B}$ are K-flat DG $k$-algebras 
(e.g.\ negatively graded and flat as $k$-modules), and 
$\tilde{A} \ar A$, $\tilde{B} \ar B$ are quasi-isomorphisms.
The ``derived category of bimodules'' should be 
$\msf{D}(\msf{DGMod} (B \otimes_{k}^{\mrm{L}} A^{\circ}))$.
It seems likely that all results in Sections 1-3 would still hold 
if we take a tilting complex to be an object of
$\msf{D}(\msf{DGMod} (B \otimes_{k}^{\mrm{L}} A^{\circ}))$, 
satisfying the appropriate conditions. However we did not check 
this.
\end{rem}

\section{Some Calculations of Tilting Complexes}

In this section we show that derived Morita equivalence reduces
to ordinary Morita equivalence when one of the algebras is local
or commutative.

\begin{lem} \label{lem2.8}
Let
$M \in \msf{D}^{-}(\cat{Mod} (B \otimes A^{\circ}))$
and
$N \in \msf{D}^{-}(\cat{Mod} (A \otimes C^{\circ}))$
for $k$-algebras $A, B, C$. Then there is a convergent K\"{u}nneth
spectral sequence
\[ \mrm{E}^{p, q}_{2} = \bigoplus_{i + j = q}
\mrm{H}^{p}(\mrm{H}^{i} M \otimes^{\mrm{L}}_{A} \mrm{H}^{j} N)
\Rightarrow \mrm{H}(M \otimes^{\mrm{L}}_{A} N) \]
in
$\cat{Mod} (B \otimes C^{\circ})$.
The filtration of each
$\mrm{H}^{n}(M \otimes^{\mrm{L}}_{A} N)$
is bounded.
If
$i_{0} \geq \sup \{ i \mid \mrm{H}^{i} M \neq 0 \}$
and
$j_{0} \geq \sup \{ j \mid \mrm{H}^{j} N \neq 0 \}$
then
\[ \mrm{H}^{i_{0}} M \otimes_{A} \mrm{H}^{j_{0}} N
\cong \mrm{H}^{i_{0} + j_{0}} (M \otimes^{\mrm{L}}_{A} N) . \]
\end{lem}

\begin{proof}
We can assume $M$ is a complex of projective
$(B \otimes A^{\circ})$-modules with
$M^{i} = 0$ for $i > i_{0}$, and similarly for $N$.
Then the usual double complex calculation applies (see \cite{ML}
Theorem XII.12.2).
In particular $\mrm{E}^{p, q}_{2} = 0$ unless
$p \leq 0$ and $q \leq i_{0} + j_{0}$.
\end{proof}

Here is a criterion for telling when we are in the classical
Morita context.

\begin{prop} \label{prop2.1}
The following conditions are equivalent for a tilting complex
$T \in \msf{D}(\cat{Mod} (B \otimes A^{\circ}))$\tup{:}
\begin{enumerate}
\rmitem{i} $T \cong P$, where $P \in \cat{Mod} (B \otimes A^{\circ})$
is invertible \tup{(}as in Theorem \tup{\ref{thm1.2})}.
\rmitem{ii} $\mrm{H}^{0} T$ is a projective $B$-module and
$\mrm{H}^{p} T = 0$ for $p \neq 0$.
\rmitem{iii} $\mrm{H}^{p} T = \mrm{H}^{p} T^{\vee} = 0$ for 
$p \neq 0$, where $T^{\vee}$ is the inverse of $T$.
\end{enumerate}
\end{prop}

\begin{proof}
(i $\Rightarrow$ ii) and (ii $\Rightarrow$ iii) are trivial.
As for (iii $\Rightarrow$ i), the shape of the K\"{u}nneth spectral
sequence shows that
\[ \begin{aligned}
\mrm{H}^{0} T \otimes_{A} \mrm{H}^{0} T^{\vee}
\cong \mrm{H}^{0} (T \otimes_{A}^{\mrm{L}} T^{\vee}) & \cong B \\[1mm]
\mrm{H}^{0} T^{\vee} \otimes_{B} \mrm{H}^{0} T
\cong \mrm{H}^{0} (T^{\vee} \otimes_{B}^{\mrm{L}} T) & \cong A .
\end{aligned} \]
\end{proof}

We call a ring $A$ {\em local} if $A / \mfrak{r}$ is a simple
artinian ring, where $\mfrak{r}$ is the Jacobson radical of $A$.
(Note that the common definition of local ring is that
$A / \mfrak{r}$ is a division ring.)

\begin{thm} \label{thm2.7}
Let $A$ and $B$ be $k$-algebras, with $A$ local, and let
$T \in \cat{D}(\cat{Mod} (B \otimes A^{\circ}))$
be a tilting complex. Then $T \cong P[n]$ for some invertible
bimodule $P$ and integer $n$.
\end{thm}

\begin{proof}
Let
$n := - \opn{max} \{ p \mid \mrm{H}^{p} T \neq 0 \}$
and
$m := - \opn{max} \{ p \mid \mrm{H}^{p} T^{\vee} \neq 0 \}$.
Then by Lemma \ref{lem1.2}, $\mrm{H}^{-n} T$ and
$\mrm{H}^{-m} T^{\vee}$ are finitely generated nonzero
modules over $A^{\circ}$ and $A$ and respectively.
Since Nakayama's Lemma holds for finitely generated $A$-modules, we 
have
$\mrm{H}^{-n} T \otimes_{A} \mrm{H}^{-m} T^{\vee} \neq 0$.
On the other hand by Lemma \ref{lem2.8} we get
\[ \mrm{H}^{-n} T \otimes_{A} \mrm{H}^{-m} T^{\vee} \cong
\mrm{H}^{-(n + m)} (T \otimes^{\mrm{L}}_{A} T^{\vee}) . \]
Since
$T \otimes^{\mrm{L}}_{A} T^{\vee} \cong B$,
we conclude that $m + n = 0$ and
\[ \mrm{H}^{-n} T \otimes_{A} \mrm{H}^{-m} T^{\vee} \cong
\mrm{H}^{0} (T \otimes^{\mrm{L}}_{A} T^{\vee}) \cong B . \]
Applying Lemma \ref{lem2.8} again we see that
\[ \mrm{H}^{-m} T^{\vee} \otimes_{B} \mrm{H}^{-n} T \cong
\mrm{H}^{0} (T^{\vee} \otimes^{\mrm{L}}_{B} T) \cong A . \]
Therefore by (ordinary) Morita equivalence it follows that
$\mrm{H}^{-n} T$ and $\mrm{H}^{-m} T^{\vee}$ are invertible bimodules.
Just as in the proof of \cite{Ye} Lemma 3.11 we find that
$\mrm{H}^{i} T = 0$ for $i \neq -n$. Taking
$P := \mrm{H}^{-n} T$ we have $T \cong P[n]$.
\end{proof}

Given a complex
$T \in \cat{D}(\cat{Mod} (B \otimes A^{\circ}))$
there are two ring homomorphisms
\begin{equation} \label{eqn2.2}
\mrm{Z}(B) \xrightarrow{\lambda_{T}}
\opn{End}_{\cat{D}(\cat{Mod} (B \otimes A^{\circ}))}(T)
\xleftarrow{\rho_{T}} \mrm{Z}(A)
\end{equation}
from the centers of $B$ and $A$, namely left and right multiplication.

\begin{prop} \label{prop2.3}
Suppose $T \in \cat{D}^{\mrm{b}}(\cat{Mod} (B \otimes A^{\circ}))$
is a tilting complex. Then the homomorphisms
$\lambda_{T}$ and $\rho_{T}$ of \tup{(\ref{eqn2.2})} are both 
bijective.
\end{prop}

\begin{proof}
Applying the functor $- \otimes_{A}^{\mrm{L}} T^{\vee}$
we get
\[ \opn{End}_{\cat{D}(\cat{Mod} (B \otimes A^{\circ}))}(T) \cong
\opn{End}_{\cat{D}(\cat{Mod} B^{\mrm{e}})}(B) \cong
\opn{End}_{\cat{Mod} B^{\mrm{e}}}(B) = \mrm{Z}(B) . \]
Since the first isomorphism sends  $\lambda_{T}$ to $\lambda_{B}$
we conclude that $\lambda_{T}$ is bijective. 
Use the functor $T^{\vee} \otimes_{B}^{\mrm{L}} -$
for $\rho_{T}$.
\end{proof}

We see that
$\mrm{Z}(A) \cong \mrm{Z}(B)$ as $k$-algebras
(cf.\ also \cite{Ri1} Proposition 9.2).

\begin{lem} \label{lem2.5}
Suppose $A$ and $B$ are $k$-algebras and
$T \in \cat{D}^{\mrm{b}}(\cat{Mod} (B \otimes A^{\circ}))$
is a tilting complex.
Let $C := \mrm{Z}(A) \cong \mrm{Z}(B)$ as in Proposition
\tup{\ref{prop2.3}},
and suppose $\tilde{C} = C S^{-1}$ is a localization of $C$ with 
respect to some multiplicative set $S \subset C$.
Define $\tilde{A} := \tilde{C} \otimes_{C} A$ and
$\tilde{B} := \tilde{C} \otimes_{C} B$.
Then 
\[ \tilde{T} := \tilde{B} \otimes_{B} T \otimes_{A} \tilde{A} \in
\cat{D}^{\mrm{b}}(\cat{Mod} (\tilde{B} \otimes 
{\tilde{A}}^{\circ})) \]
is a tilting complex, with inverse
$\tilde{T}^{\vee} := 
\tilde{A} \otimes_{A} T^{\vee} \otimes_{B} \tilde{B}$.
\end{lem}

\begin{proof}
By Proposition \ref{prop2.3} the cohomology bimodules 
$\mrm{H}^{p} T$ are all $C$-central 
(even though $T$ itself need not be $C$-central!).
From the flatness of $A \ar \tilde{A}$ and $B \ar \tilde{B}$,
and using the fact that 
$\tilde{C} \otimes_{C} \tilde{C} \cong \tilde{C}$,
we conclude that
\[ \mrm{H}^{p} \tilde{T} \cong \tilde{B} \otimes_{B} \mrm{H}^{p} T
\otimes_{A} \tilde{A} \cong
\tilde{C} \otimes_{C} \mrm{H}^{p} T . \]
Hence
$\tilde{B} \otimes_{B} T \ar \tilde{T}$
and
$T \otimes_{A} \tilde{A} \ar \tilde{T}$
are isomorphisms in
$\msf{D}(\cat{Mod}(B \otimes A^{\circ}))$.

The functor $\mrm{R}\opn{Hom}_{B}(-, B)$
gives rise to an isomorphism
\[ \opn{End}_{\msf{D}(\cat{Mod}(B \otimes A^{\circ}))}(T)^{\circ}
\cong
\opn{End}_{\msf{D}(\cat{Mod}(A \otimes B^{\circ}))}(T^{\vee}) \]
which exchanges $\rho$ and $\lambda$. Therefore the
$\mrm{H}^{p} T^{\vee}$ are also $C$-central, and as above
$\tilde{A} \otimes_{A} T^{\vee} \cong \tilde{T}^{\vee} \cong
T^{\vee} \otimes_{B} \tilde{A}$.
We see that
\[ \begin{aligned}
\tilde{T}^{\vee} \otimes^{\mrm{L}}_{\tilde{B}}
\tilde{T} & \cong
(\tilde{A} \otimes_{A} T^{\vee} \otimes_{B} \tilde{B})
\otimes^{\mrm{L}}_{\tilde{B}}
(\tilde{B} \otimes_{B} T \otimes_{A} \tilde{A}) \\
& \cong \tilde{A} \otimes_{A} ({T^{\vee}} \otimes^{\mrm{L}}_{B} T)
\otimes_{A} \tilde{A} \cong \tilde{A}
\end{aligned} \]
and likewise
$\tilde{T} \otimes^{\mrm{L}}_{\tilde{A}}
\tilde{T}^{\vee} \cong \tilde{B}$.
\end{proof}

In Morita equivalence
(i.e.\ Theorem \ref{thm1.2}), if the ring $A$ is commutative
then the isomorphism $A \cong \mrm{Z}(B)$ makes the
invertible bimodule $P$ $A$-central. Since $B \cong \opn{End}_{A}(P)$
it is an Azumaya algebra over $A$.
The next theorem says that in the commutative case, derived
Morita equivalence gives nothing new.

\begin{thm} \label{thm2.4}
Let $A$ and $B$ be $k$-algebras, with $A$ commutative.
If $A$ and $B$ are derived Morita equivalent, then they are
Morita equivalent.
\end{thm}

\begin{proof}
By Corollary \ref{cor1.6} there exists a tilting complex
$T \in \cat{D}^{\mrm{b}}(\cat{Mod} (B \otimes A^{\circ}))$.
Let $T^{\vee}$ be its inverse.
Write $A = A_{1} \times \cdots \times A_{m}$
with $\opn{Spec} A_{i}$ connected. Let $B_{i} := A_{i} \otimes_{A} B$.
Then by Lemma \ref{lem2.5},
$B_{i} \otimes_{B} T \otimes_{A} A_{i}$ is a tilting complex
in $\cat{D}^{\mrm{b}}(\cat{Mod} (B_{i} \otimes A_{i}^{\circ}))$.
Thus we may assume $\opn{Spec} A$ is connected.

Pick a prime ideal
$\mfrak{p} \in \opn{Spec} A$, let $A_{\mfrak{p}}$ be the local ring,
$B_{\mfrak{p}} := A_{\mfrak{p}} \otimes_{A} B$
and
$T_{\mfrak{p}} := B_{\mfrak{p}} \otimes_{B} 
T \otimes_{A} A_{\mfrak{p}}$.
By Lemma \ref{lem2.5} the complex $T_{\mfrak{p}}$
is a tilting complex in
$\cat{D}^{\mrm{b}}(\cat{Mod} (B_{\mfrak{p}} \otimes
A_{\mfrak{p}}^{\circ}))$,
with inverse 
$T^{\vee}_{\mfrak{p}} := B_{\mfrak{p}} \otimes_{B} 
T^{\vee} \otimes_{A} A_{\mfrak{p}}$.
Define integers $n(\mfrak{p})$ and $m(\mfrak{p})$ by
$n(\mfrak{p}) := - \opn{max}
\{ i \mid \mrm{H}^{i} T_{\mfrak{p}} \neq 0 \}$
and
$m(\mfrak{p}) := - \opn{max}
\{ i \mid \mrm{H}^{i} T^{\vee}_{\mfrak{p}} \neq 0 \}$.
As in the proof of Theorem \ref{thm2.7},
$\mrm{H}^{-n(\mfrak{p})} T_{\mfrak{p}} \cong
A_{\mfrak{p}} \otimes_{A} \mrm{H}^{-n(\mfrak{p})} T$
is an invertible $B_{\mfrak{p}}$-$A_{\mfrak{p}}$-bimodule, 
$\mrm{H}^{i} T_{\mfrak{p}} \cong
A_{\mfrak{p}} \otimes_{A} \mrm{H}^{i} T = 0$
for $i \neq -n(\mfrak{p})$, and $m(\mfrak{p}) + n(\mfrak{p}) = 0$.

Next consider prime ideals $\mfrak{p} \subset \mfrak{q}$.
The previous paragraph implies that
$A_{\mfrak{p}} \otimes_{A} \mrm{H}^{i} T = 0$
for $i \neq -n(\mfrak{q})$, and hence
$n(\mfrak{q}) = n(\mfrak{p})$. Because $\opn{Spec} A$ is connected
we conclude that $n(\mfrak{p}) = n$ is constant, and so
$\mrm{H}^{i} T = 0$ for $i \neq -n$.
Likewise $m(\mfrak{p}) = m = -n$ and 
$\mrm{H}^{i} T^{\vee} = 0$ for $i \neq -m$. 
By Proposition \ref{prop2.1} we see that the $A$-central bimodule
$P := \mrm{H}^{-n} T$ is invertible.
\end{proof}

Here is a corollary to Theorem \ref{thm2.7}:

\begin{cor} \label{cor2.7}
Let $A$ and $B$ be $k$-algebras, and
$A \cong A_{1} \times \cdots \times A_{m}$
with $A_{i}$ local.
If $A$ and $B$ are derived Morita equivalent, then they are
Morita equivalent.
\end{cor}

\begin{proof}
Let $C := \mrm{Z}(A)$, so $C = C_{1} \times \cdots \times C_{m}$.
By Lemma \ref{lem2.5} every
$B_{i} := C_{i} \otimes_{C} B$ is derived Morita equivalent to $A_{i}$.
Now use  Theorem \ref{thm2.7}.
\end{proof}

\begin{rem} \label{rem2.1}
R. Rouquier and A. Zimmermann have independently obtained similar 
results to our Theorems \ref{thm2.7} and \ref{thm2.4}, but only in 
a special case: when $A$ and $B$ are orders over a Dedekind domain 
$k$. They also considered the derived Picard group, which they 
denoted by $\mrm{TrPic}(A)$. See \cite{Zi}. 
\end{rem}

\section{The Derived Picard Group}

Let us concentrate now on the case $A = B$.
Recall that the $k$-central noncommutative {\em Picard group}
of $A$ is
\[ \opn{Pic}(A) = \opn{Pic}_{k}(A) :=
\frac{ \{ \text{invertible bimodules }
L \in \cat{Mod} A^{\mrm{e}} \} }{
\text{isomorphism} } . \]
According to Corollary \ref{cor1.1} the next definition makes sense:

\begin{dfn} \label{dfn3.2}
Define the {\em derived Picard group} of $A$ \tup{(}relative to
$k$\tup{)} to be
\[ \opn{DPic}(A) = \opn{DPic}_{k}(A) :=
\frac{ \{ \text{tilting complexes}\
T \in \cat{D}^{\mrm{b}}(\cat{Mod} A^{\mrm{e}}) \} }{
\text{isomorphism} } \]
with identity element $A$, product
$(T_{1}, T_{2}) \mapsto T_{1} \otimes_{A}^{\mrm{L}} T_{2}$
and inverse
$T \mapsto T^{\vee}$.
\end{dfn}

The group $\opn{DPic}(A)$ contains a copy of $\mbb{Z}$ in its center,
as $n \mapsto A[n]$.
$\opn{DPic}(A)$ also contains a subgroup isomorphic to $\opn{Pic}(A)$,
which is characterized in Proposition \ref{prop2.1}.
Note that both $\opn{Pic}(A)$ and $\opn{DPic}(A)$ depend on $k$.

\begin{rem}
If $A$ is commutative we denote by $\opn{Pic}_{A}(A)$ the usual
commutative Picard group, namely the isomorphism classes of central
projective modules of rank $1$. It is a subgroup of
$\opn{Pic}(A)$ (cf.\ Proposition \ref{prop3.3}).
\end{rem}

Let us now state some facts about invertible bimodules
(which are probably well known, but we found no references).
Denote by $\opn{Aut}(A)$ the group of $k$-algebra automorphisms of
$A$. For $\sigma \in \opn{Aut}(A)$ let $A_{\sigma}$ be
the bimodule which is free over $A$ and $A^{\circ}$ with basis $e$,
and $e \cdot a = \sigma(a) \cdot e$, $a \in A$.

\begin{lem} \label{lem3.1}
\begin{enumerate}
\item $\sigma \mapsto A_{\sigma}$ is a group homomorphism
$\opn{Aut}(A) \ar \opn{Pic}(A)$
with kernel the subgroup $\opn{Inn}(A)$ of inner automorphisms.
\item Suppose $L$ is an invertible $A$-bimodule which is free of rank
$1$ as left module. Then $L \cong A_{\sigma}$ as bimodules for some
$\sigma \in \opn{Aut}(A)$.
\item If $A$ is local then any
invertible bimodule $L$ is free of rank $1$ over $A$.
\end{enumerate}
\end{lem}

\begin{proof}
1.\ A bimodule isomorphism
$A_{\tau} \iso A_{\sigma}$  sends the basis $e$ of
$A_{\tau}$ to $u \cdot e \in A_{\sigma}$, where $u$ is
a unit of $A$, and conjugation by $u$ is $\tau \sigma^{-1}$.
\newline
2.\ Choose an $A$-basis $e$ of $L$. Then $\phi \mapsto \phi(e)$
is a bijection $\opn{End}_{A}(L) \iso L$. But
since right multiplication induces an isomorphism
$A^{\circ} \iso \opn{End}_{A}(L)$ this shows that $e$
is also a basis of $L$ as $A^{\circ}$-module. Conjugation by $e$ is
$\sigma$.
\newline
3.\ By Nakayama's Lemma it is enough to prove that
$W := K \otimes_{A} L$ is free of rank $1$ over
$K := A / \mfrak{r}$. First one checks that
$W \cong K \otimes_{A} L \otimes_{A} K \cong L \otimes_{A} K$
as $K^{\mrm{e}}$-modules. Hence $W$ is an invertible bimodule over $K$.
Since $K \cong \mrm{M}_{n}(D)$
for a division algebra $D$, by Morita equivalence
$W \cong \mrm{M}_{n}(V)$ as $K^{\mrm{e}}$-modules, where $V$ is
an invertible bimodule over $D$. It remains to prove that the free
$D$-module $V$ has rank $1$. But if $V \cong D^{l}$ as left modules,
then
$D \cong V^{\vee} \otimes_{D} V \cong (V^{\vee})^{l}$,
so $l = 1$.
\end{proof}

The next propositions analyze $\opn{DPic}(A)$ in the semilocal
and in the commutative cases.

\begin{prop} \label{prop3.4}
Suppose $A \cong A_{1} \times \cdots \times A_{m}$ where every
$A_{i}$ is a local $k$-algebra.
Then
\[ \begin{aligned}
\opn{DPic}_{k}(A) & \cong \mbb{Z}^{m} \times \opn{Pic}_{k}(A) \\[1mm]
\opn{Pic}_{k}(A) & \cong \opn{Out}_{k}(A) .
\end{aligned} \]
\end{prop}

\begin{proof}
This is an immediate consequence of Corollary \ref{cor2.7},
Theorem \ref{thm2.7} and and Lemma \ref{lem3.1}.
\end{proof}

\begin{prop} \label{prop3.3}
Suppose $A$ is a commutative ring. Then
\[ \begin{aligned}
\opn{DPic}_{k}(A) & \cong \mbb{Z}^{m} \times \opn{Pic}_{k}(A) \\[1mm]
\opn{Pic}_{k}(A) & \cong \opn{Aut}_{k}(A) \ltimes \opn{Pic}_{A}(A) ,
\end{aligned} \]
where $m$ is the number of connected components of
$\opn{Spec} A$.
\end{prop}

\begin{proof}
Let $A \cong A_{1} \times \cdots A_{m}$ be the decomposition
of $A$ according to the connected components of $\opn{Spec} A$.
Given a tilting complex
$T \in \cat{D}^{\mrm{b}}(\cat{Mod} A^{\mrm{e}})$,
Theorem \ref{thm2.4} says that
$T \cong P \otimes_{A} S$, where $P$ is an invertible bimodule and
$S = A_{1}[n_{1}] \times \cdots \times A_{m}[n_{m}]$
with $n_{i} \in \mbb{Z}$. So
$\opn{DPic}_{k}(A) \cong \mbb{Z}^{m} \times \opn{Pic}_{k}(A)$.
Next let $\sigma \in \opn{Aut}(A)$ be the automorphism determined
by $P$ (cf.\ Proposition \ref{prop2.3}). Then
$L := P \otimes_{A} A_{\sigma^{-1}}$
is a central invertible bimodule over $A$.
\end{proof}

Assume $A$ is noetherian. Let
$\mrm{K}_{0}(A) = \mrm{K}_{0}(\cat{Mod}_{\mrm{f}} A)$
be the Grothendieck group of $A$. For any
$M \in \cat{Mod}_{\mrm{f}} A$ let $[M]$ be its class in
$\mrm{K}_{0}(A)$. Then
$M \mapsto [M] := \sum (-1)^{p} [\mrm{H}^{p}M]$
is a well defined function
$\msf{D}^{\mrm{b}}_{\mrm{f}}(\cat{Mod} A) \ar \mrm{K}_{0}(A)$.
Since $\opn{DPic}(A)$ acts on
$\msf{D}^{\mrm{b}}_{\mrm{f}}(\cat{Mod} A)$
by auto-equivalences it acts also on $\mrm{K}_{0}(A)$.
Let $s \in \opn{DPic}(A)$ be the class of $A[1]$, which acts on
$\msf{D}^{\mrm{b}}_{\mrm{f}}(\cat{Mod} A)$ by a shift in degree.
Then:

\begin{prop} \label{prop3.5}
There is a canonical group homomorphism
\[ \chi_{0} : \opn{DPic}(A) \ar \opn{Aut}_{\mbb{Z}}(\mrm{K}_{0}(A)) \]
with $\chi_{0}(s) = -1$.
\end{prop}

Actually there are two more objects one can associate to $A$
which are related to the representation $\chi_{0}$.

The first is the noncommutative Grothendieck ring
$\opn{K}^{0}(A) = \opn{K}^{0}_{k}(A)$,
which is a rather obvious generalization of the commutative
$\opn{K}^{0}(A)$.
Let $X$ be the set of isomorphism classes of $A^{\mrm{e}}$-modules
$T$ which are finitely generated projective on both sides.
Define $F$ to be the free abelian group with basis $X$.
As abelian group, $\opn{K}^{0}(A)$ is the quotient of $F$
by the subgroup generated by the elements
$[T_{0}] - [T_{1}] + [T_{2}]$,
for every short exact sequence
$0 \ar T_{0} \ar T_{1} \ar T_{2} \ar 0$ in $\cat{Mod} A^{\mrm{e}}$
with $T_{i} \in X$.
Multiplication is
$[T_{1}] \cdot [T_{2}] := [T_{1} \otimes_{A} T_{2}]$,
and $1$ is $[A]$. The Grothendieck group $\opn{K}_{0}(A)$
is a left module over $\opn{K}^{0}(A)$, by
$[T] \cdot [M] := [T \otimes_{A} M]$
for $M \in \cat{Mod} A$, and there is a group homomorphism
$\opn{Pic}(A) \ar \opn{K}^{0}(A)^{\times}$.

All the above works for complexes too: take $X$ to be the set of
isomorphism classes of complexes
$T \in \msf{D}^{\mrm{b}}(\cat{Mod} A^{\mrm{e}})$
which are perfect on both sides. Define $F$ as before, and let
$\opn{DK}^{0}(A) = \opn{DK}^{0}_{k}(A)$
be the quotient of $F$
by the subgroup generated by the elements
$[T_{0}] - [T_{1}] + [T_{2}]$,
for every triangle
$T_{0} \ar T_{1} \ar T_{2} \xrightarrow{+1}$ in
$\msf{D}^{\mrm{b}}(\cat{Mod} A^{\mrm{e}})$
with $T_{i} \in X$. Multiplication is
$[T_{1}] \cdot [T_{2}] := [T_{1} \otimes^{\mrm{L}}_{A} T_{2}]$.
The Grothendieck group $\opn{K}_{0}(A)$
is a left module over $\opn{DK}^{0}(A)$, by
$[T] \cdot [M] := [T \otimes^{\mrm{L}}_{A} M]$
for $M \in \msf{D}^{\mrm{b}}_{\mrm{f}}( \cat{Mod} A)$.
There is a ring homomorphism
$\opn{K}^{0}(A) \ar \opn{DK}^{0}(A)$,
and a group homomorphism
$\opn{DPic}(A) \ar \opn{DK}^{0}(A)^{\times}$.
To summarize:

\begin{prop} \label{prop3.8}
$\opn{DK}^{0}(A)$ is a ring and
$\opn{K}_{0}(A)$ is a left $\opn{DK}^{0}(A)$-module.
There is a group homomorphism
$\chi^{0} : \opn{DPic}(A) \ar \opn{DK}^{0}(A)^{\times}$
with $\chi^{0}(s) = -1$,
and $\chi_{0}$ factors through $\chi^{0}$.
\end{prop}

\begin{rem}
We did not analyze the dependence of various objects,
such as the group $\opn{DPic}(A)$, on the base field $k$.
\end{rem}

\section{Classification of Dualizing Complexes}

In this section we assume that $A$ is a (left and right) noetherian
$k$-algebra.
Dualizing complexes over commutative rings were introduced
in \cite{RD}. The noncommutative version below first appeared in
\cite{Ye} (where connected graded algebras were considered).

\begin{dfn} \label{dfn4.1}
A complex $R \in \cat{D}^{\mrm{b}}(\cat{Mod} A^{\mrm{e}})$ is
called {\em dualizing} if
\begin{enumerate}
\rmitem{i} $R$ has finite injective dimension over $A$ and over
$A^{\circ}$.
\rmitem{ii} $R$ has finitely generated cohomology modules
over $A$ and over $A^{\circ}$.
\rmitem{iii} Then canonical morphisms
$A \ar \mrm{R} \opn{Hom}_{A}(R, R)$
and
$A \ar \mrm{R} \opn{Hom}_{A^{\circ}}(R, R)$
in $\cat{D}^{\mrm{b}}(\cat{Mod} A^{\mrm{e}})$ are isomorphisms.
\end{enumerate}
\end{dfn}

Condition (i) is equivalent to having an isomorphism
$R \cong I \in \cat{D}^{\mrm{b}}(\cat{Mod} A^{\mrm{e}})$,
where $I$ is a bounded complex and each $I^{q}$ is injective
over $A$ and over $A^{\circ}$.
Note that this definition is left-right symmetric
(i.e.\ remains equivalent after exchanging $A$ and $A^{\circ}$).

Given a dualizing complex $R$ the associated duality functors
are
\[ \begin{aligned}
D := \mrm{R} \opn{Hom}_{A}(-, R) : &
\cat{D}(\cat{Mod} A)^{\circ} \ar \cat{D}(\cat{Mod} A^{\circ}) \\
D^{\circ} := \mrm{R} \opn{Hom}_{A^{\circ}}(-, R) : &
\cat{D}(\cat{Mod} A^{\circ})^{\circ} \ar \cat{D}(\cat{Mod} A) .
\end{aligned} \]

For a $k$-algebra $B$ let
$\cat{D}_{(\mrm{f}, \, )}(\cat{Mod} (A \otimes B))$
denote the full subcategory of \linebreak
$\cat{D}(\cat{Mod} (A \otimes B))$
whose objects are the complexes $M$ s.t.\ for all $q$,
$\mrm{H}^{q} M$ is a finitely generated $A$-module.
Likewise define
$\cat{D}_{(\, , \mrm{f})}(\cat{Mod} (A \otimes B))$
and
$\cat{D}_{(\mrm{f}, \mrm{f})}(\cat{Mod} (A \otimes B))$.
Thus condition (ii) of Definition \ref{dfn4.1} says that
$R \in \cat{D}_{(\mrm{f}, \mrm{f})}(\cat{Mod} A^{\mrm{e}})$.

\begin{prop} \label{prop4.1}
Let
$R \in \cat{D}^{\mrm{b}}(\cat{Mod} A^{\mrm{e}})$
be a dualizing complex, and let $B$ be any $k$-algebra.
\begin{enumerate}
\item For any
$M \in \cat{D}_{(\mrm{f}, \, )}(\cat{Mod} (A \otimes B^{\circ}))$
one has
$D M \in \cat{D}_{(\, , \mrm{f})}(\cat{Mod} (B \otimes A^{\circ}))$,
and there is a functorial isomorphism
$M \cong D^{\circ} D M$.
Therefore $D$ and $D^{\circ}$ determine an equivalence
\[ \cat{D}_{(\mrm{f}, \, )}(\cat{Mod} (A \otimes B^{\circ}))^{\circ}
\longleftrightarrow
\cat{D}_{(\, , \mrm{f})}(\cat{Mod} (B \otimes A^{\circ})) . \]
\item Let
$M \in
\cat{D}^{-}_{(\mrm{f}, \, )}(\cat{Mod} (A \otimes B^{\circ}))$
and
$N \in
\cat{D}_{(\mrm{f}, \, )}(\cat{Mod} (A \otimes B^{\circ}))$.
Then there is a bifunctorial isomorphism
\[ \mrm{R} \opn{Hom}_{A}(M, N) \cong
\mrm{R} \opn{Hom}_{A^{\circ}}(D N, D M)  \]
in $\cat{D}(\cat{Mod} B)$.
\end{enumerate}
\end{prop}

\begin{proof}
1.\  This is slightly stronger than \cite{Ye} Lemma 3.5.
By adjunction we get a functorial morphism
$M \ar D^{\circ} D M$ in
$\cat{D}(\cat{Mod} (A \otimes B))$. Now we can forget $B$.
Since the functors $D$ and $D^{\circ} D$ are way-out in both
directions,
$D A = R \in \cat{D}_{\mrm{f}}(\cat{Mod} A^{\circ})$
and $D^{\circ} D A \cong A$,
the claim follows from \cite{RD} Propositions I.7.1 and I.7.3
and their opposite forms.
\newline
2.\ We can assume $M$ is a bounded above complex of projective
$(A \otimes B^{\circ})$-modules and $R$ is a bounded below
complex of injective $A^{\mrm{e}}$-modules. Since
$\opn{Hom}_{A}(M, R)$ is a bounded below complex of injective
$A^{\circ}$-modules, we get a morphism
\[ \begin{split}
\mrm{R} \opn{Hom}_{A}(M, N) & = \opn{Hom}_{A}(M, N) \ar \\
& \opn{Hom}_{A^{\circ}}(\opn{Hom}_{A}(N, R), \opn{Hom}_{A}(M, R)) =
\mrm{R} \opn{Hom}_{A^{\circ}}(D N, D M)
\end{split} \]
in $\cat{D}(\cat{Mod} B)$, which is functorial in $M, N$. In order
to prove it is an isomorphism we can forget $B$. Applying
$\mrm{H}^{q}$ we get a homomorphism
$\opn{Hom}_{\cat{D}(\cat{Mod} A)}(M, N[q]) \ar
\opn{Hom}_{\cat{D}(\cat{Mod} A^{\circ})}(D(N[q]), D M)$,
which by part 1 is bijective.
\end{proof}

\begin{exa}
If $A$ is a Gorenstein algebra, that is the bimodule $A$ has
finite injective dimension over $A$ and over $A^{\circ}$, then
$R = A$ is a dualizing complex.
\end{exa}

\begin{rem}
One can weaken the noetherian assumption. If $A$ is a
coherent ring then the category of coherent (i.e.\ finitely presented)
$A$-modules is abelian, so we can work with
$\cat{D}^{\mrm{b}}_{\mrm{c}}(\cat{Mod} A)$ etc.
Perhaps a reasonable theory can be developed for any algebra $A$
if one works with Illusie's pseudo-coherent complexes
(cf.\ \cite{SGA6} Expos\'{e}e I).
\end{rem}

Let
$\cat{D}^{-}(\cat{Mod} A)_{\mrm{fpd}}$
(resp.\ $\cat{D}^{-}(\cat{Mod} A)_{\mrm{fTd}}$,
$\cat{D}^{+}(\cat{Mod} A)_{\mrm{fid}}$)
be the category of complexes with finite projective (resp.\ Tor,
injective) dimension. Since $A$ is noetherian, we have
\[ \cat{D}^{\mrm{b}}_{\mrm{f}}(\cat{Mod} A)_{\mrm{fTd}} =
\cat{D}^{\mrm{b}}_{\mrm{f}}(\cat{Mod} A)_{\mrm{fpd}} =
\cat{D}^{\mrm{b}}(\cat{Mod} A)_{\mrm{perf}}
\subset \cat{D}^{\mrm{b}}(\cat{Mod} A) . \]

\begin{thm} \label{thm4.1}
\begin{enumerate}
\item Suppose $R_{1}$ is a dualizing complex and $T$ is a tilting
complex. Then $R_{2} := R_{1} \otimes^{\mrm{L}}_{A} T$
is dualizing, and
$T \cong \mrm{R} \opn{Hom}_{A}(R_{1}, R_{2})$.
\item Conversely, suppose $R_{1}$  and $R_{2}$ are dualizing complexes.
Then
$T := $ \newline
$\mrm{R} \opn{Hom}_{A}(R_{1}, R_{2})$
is a tilting complex $T$, and
$R_{2} \cong R_{1} \otimes^{\mrm{L}}_{A} T$.
\item Let $R$ be a dualizing complex. Then the associated duality
functors $D$ and $D^{\circ}$
induce an equivalence
\[ \cat{D}^{\mrm{b}}_{\mrm{f}}(\cat{Mod} A)_{\mrm{fpd}}^{\circ}
\longleftrightarrow
\cat{D}^{\mrm{b}}_{\mrm{f}}(\cat{Mod} A^{\circ})_{\mrm{fid}} . \]
\end{enumerate}
\end{thm}

\begin{proof}
1.\
Clearly $R_{2}$ is bounded.
Next let us prove that each $\mrm{H}^{n} R_{2}$ is a
finitely generated module over
$A$. Consider the K\"{u}nneth spectral sequence
\[ E^{pq}_{2} = \bigoplus_{i + j = q}
\mrm{H}^{p}(\mrm{H}^{i} R_{1} \otimes^{\mrm{L}}_{A}
\mrm{H}^{j} T) \Rightarrow
\mrm{H}(R_{1} \otimes^{\mrm{L}}_{A} T)
= \mrm{H} R_{2} . \]
Using a resolution of $\mrm{H}^{j} T$ by finitely generated
flat $A$-modules one easily sees that
$\mrm{H}^{p}(\mrm{H}^{i} R_{1} \otimes^{\mrm{L}}_{A}
\mrm{H}^{j} T)$
is finitely generated over $A$.
Since the filtration on $\mrm{H}^{n} R_{2}$ is bounded it follows
that this too is a finitely generated $A$-module.
Finiteness over $A^{\circ}$ is proved similarly.

Given
$M \in \cat{D}^{\mrm{b}}(\cat{Mod} A)$
there is a natural isomorphism
\[ \mrm{R} \opn{Hom}_{A}(M, R_{1} \otimes^{\mrm{L}}_{A} T) \cong
\mrm{R} \opn{Hom}_{A}(M, R_{1}) \otimes^{\mrm{L}}_{A} T . \]
If
$M \in \cat{D}^{\mrm{b}}(\cat{Mod} A^{\mrm{e}})$
there is also a natural isomorphism
\[ \mrm{R} \opn{Hom}_{A}(M, R_{1} \otimes^{\mrm{L}}_{A} T) \cong
\mrm{R} \opn{Hom}_{A}(M \otimes^{\mrm{L}}_{A} T^{\vee} , R_{1}) \]
where
$T^{\vee} := \mrm{R} \opn{Hom}_{A}(T, A)$.
Therefore $R_{2}$ has finite injective dimension over $A$,
$A \cong \mrm{R} \opn{Hom}_{A}(R_{2}, R_{2})$
and
$T \cong \mrm{R} \opn{Hom}_{A}(R_{1}, R_{2})$.
There is also a natural isomorphism
\[ \mrm{R} \opn{Hom}_{A^{\circ}}(N, R_{1} \otimes^{\mrm{L}}_{A} T) 
\cong \mrm{R} \opn{Hom}_{A^{\circ}}(N \otimes^{\mrm{L}}_{A}
T^{\vee}, R_{1})  \]
for
$N \in \cat{D}^{\mrm{b}}(\cat{Mod} A^{\circ})$,
so $R_{2}$ has finite injective dimension over $A^{\circ}$
and
$A \cong \mrm{R} \opn{Hom}_{A^{\circ}}(R_{2}, R_{2})$.

\medskip \noindent
2.\ By the proof of \cite{Ye} Theorem 3.9, $T$ is a tilting
complex, and by ibid.\ Lemma 3.10,
$R_{2} \cong R_{1} \otimes^{\mrm{L}}_{A} T$.

\medskip \noindent
3.\ Just like \cite{RD} Proposition IV.2.6.
\end{proof}

The theorem says that 
$(R, T) \mapsto R \otimes^{\mrm{L}}_{A} T$
is a right action of $\opn{DPic}(A)$ on the set of isomorphism classes
of dualizing complexes. By symmetry there is a left action 
$(T, R) \mapsto T \otimes^{\mrm{L}}_{A} R$. 
As a corollary we get the classification of isomorphism classes
of dualizing complexes. 

\begin{cor} \label{cor4.1}
If the set
\[ \frac{ \{ \text{\tup{dualizing complexes}}\
R \in \cat{D}^{\mrm{b}}(\cat{Mod} A^{\mrm{e}}) \} }{
\text{\tup{isomorphism}} }  \]
is nonempty, then the left and right actions of the group 
$\opn{DPic}(A)$ on it are transitive with trivial stabilizers. 
\end{cor}

\begin{prob}
In Propositions \ref{prop3.3} and \ref{prop3.4} we have seen that when
$A$ is commutative or local, the group $\opn{DPic}(A)$ consists of
familiar ingredients - $\opn{Pic}_{A}(A)$, $\opn{Aut}(A)$ and the
trivial copy of $\mbb{Z}$ (cf.\ also Section 6).
On the other hand $\opn{DPic}(A)$ classifies the isomorphism classes
of dualizing complexes. Now it is known in commutative algebraic
geometry that dualizing complexes are in close relation to
localization. For instance, a ring with a dualizing complex is
catenary; a dualizing complex can be represented by a residual complex,
which is a sum of local cohomology modules (see \cite{RD}).
This leads us to ask whether some obstructions to localization
can be found in $\opn{DPic}(A)$ when $A$ is noncommutative?
More specifically, is there a relation between the group structure
of $\opn{DPic}(A)$ and the link graph of maximal ideals in
$\opn{Spec} A$?
\end{prob}

\section{Rigid Dualizing Complexes}

In this section we use the action of the group $\opn{DPic}(A)$
on the set of isomorphism classes of dualizing complexes to study
certain properties of dualizing complexes. In particular we shall
be interested in rigid dualizing complexes, which were recently
introduced by M.\ Van den Bergh.
As in Section 4, $A$ is a noetherian $k$-algebra.

First we need some notational conventions on modules with
multiple actions.
For an element $a \in A$ we denote by $a^{\circ} \in A^{\circ}$
the same element. Thus for $a_{1}, a_{2} \in A$,
$a_{1}^{\circ} \cdot a_{2}^{\circ} =
(a_{2} \cdot a_{1})^{\circ} \in A^{\circ}$.
With this notation if $M$ is a right $A$-module then the left
$A^{\circ}$ action is
$a^{\circ} \cdot m = m \cdot a$, $m \in M$.
The algebra $A^{\mrm{e}}$ has an involution
$A^{\mrm{e}} \ar (A^{\mrm{e}})^{\circ}$,
$a_{1} \otimes a_{2}^{\circ} \mapsto a_{2} \otimes a_{1}^{\circ}$
which allows us to regard every left $A^{\mrm{e}}$-module $M$
as a right $A^{\mrm{e}}$-module in a consistent way:
\[ (a_{1} \otimes a_{2}^{\circ}) \cdot m =
(a_{2} \otimes a_{1}^{\circ})^{\circ} \cdot m =
m \cdot (a_{2} \otimes a_{1}^{\circ}) = a_{1} \cdot m \cdot a_{2} . \]

Given an $(A \otimes B)$-module $M$ and a $(B \otimes A)$-module $N$
we define a mixed action of
$A^{\mrm{e}} \otimes B^{\mrm{e}}$
on the tensor product $M \otimes N$ as follows.
$A^{\mrm{e}}$ acts on $M \otimes N$  by the outside action
\[ (a_{1} \otimes a_{2}^{\circ}) \cdot (m \otimes n) :=
(a_{1} \cdot m) \otimes (n \cdot a_{2}) , \]
whereas $B^{\mrm{e}}$ acts on $M \otimes N$  by the inside action
\[ (b_{1} \otimes b_{2}^{\circ}) \cdot (m \otimes n) :=
(m \cdot b_{2}) \otimes (b_{1} \cdot n) . \]
By default we regard the outside action as a left action and the
inside action as a right action. If $A = B$ and $M = N$ then
the two actions by $A^{\mrm{e}}$ on $M \otimes M$ are interchanged
the involution
$m_{1} \otimes m_{2} \mapsto m_{2} \otimes m_{1}$. However for
the sake of definiteness in this case, given an $A^{\mrm{e}}$-module
$L$,
$\opn{Hom}_{A^{\mrm{e}}}(L, M \otimes M)$ shall refer to
homomorphisms $L \ar M \otimes M$ which are $A^{\mrm{e}}$-linear
w.r.t.\ the outside action.

The next definition is due to M.\ Van den Bergh in \cite{VdB}.

\begin{dfn}
A {\em rigid dualizing complex} over $A$ is a pair $(R, \rho)$ where
$R$ is a dualizing complex and
\[ \rho : R \iso
\mrm{R} \opn{Hom}_{A^{\mrm{e}}}(A, R \otimes R) \]
is an isomorphism in $\cat{D}(\cat{Mod}(A^{\mrm{e}}))$.
\end{dfn}

Van den Bergh proved that any two
rigid dualizing complexes are isomorphic. We improve slightly:

\begin{thm} \label{thm5.1}
Suppose $(R_{1}, \rho_{1})$ and $(R_{2}, \rho_{2})$
are two rigid dualizing complexes. Then there is a {\em unique}
isomorphism
$\phi : R_{1} \iso R_{2}$
in $\cat{D}(\cat{Mod}(A^{\mrm{e}}))$ making the diagram
\begin{equation} \label{eqn1}
\begin{CD}
R_{1} @> \rho_{1} >>
\mrm{R} \opn{Hom}_{A^{\mrm{e}}}
(A, R_{1} \otimes R_{1}) \\
@V \phi VV @V \phi \otimes \phi VV \\
R_{2} @> \rho_{2} >>
\mrm{R} \opn{Hom}_{A^{\mrm{e}}}
(A, R_{2} \otimes R_{2})
\end{CD}
\end{equation}
commute.
\end{thm}

First we need:

\begin{lem} \label{lem5.2}
Let $R$ be a dualizing complex. Then the two ring homomorphisms
$\lambda_{R}, \rho_{R} :
\mrm{Z}(A) \ar \opn{End}_{\cat{D}(\cat{Mod} A^{\mrm{e}})}(R)$,
namely left and right multiplication, are bijective.
\end{lem}

\begin{proof}
The proof is similar to Proposition \ref{prop2.3}.
Define functors
$D := \mrm{R} \opn{Hom}_{A}(-, R)$
and
$D^{\circ} := \mrm{R} \opn{Hom}_{A^{\circ}}(-, R)$.
Since
$A \cong D^{\circ} D A \cong D^{\circ} R$
in $\cat{D}(\cat{Mod} A^{\mrm{e}})$ it follows (by applying 
$D^{\circ}$) that
\[ \opn{Hom}_{\cat{D}(\cat{Mod} A^{\mrm{e}})}(R, R) \cong
\opn{Hom}_{\cat{D}(\cat{Mod} A^{\mrm{e}})}(A, A)^{\circ} . \]
This sends the left action $\lambda_{R}$ of $\mrm{Z}(A)$ on $R$
to the right action $\rho_{A}$ of $\mrm{Z}(A)$ on $A$. But
\[ \opn{End}_{\cat{D}(\cat{Mod} A^{\mrm{e}})}(A) =
\opn{End}_{\cat{Mod} A^{\mrm{e}}}(A) = \mrm{Z}(A)  \]
(via $\lambda_{A} = \rho_{A}$). Hence $\lambda_{R}$ is bijective.
Do the same for $\rho_{R}$.
\end{proof}

\begin{proof}[Proof of the theorem]
Suppose we are given some isomorphism
$\phi' : R_{1} \iso R_{2}$. Let
$\psi \in \opn{Aut}(R_{2})$ satisfy
\[ (\phi' \otimes \phi') \rho_{1} = \rho_{2} \psi \phi' , \]
and define $\phi := \psi^{-1} \phi'$.
By Lemma \ref{lem5.2} there are elements
$a, b \in \mrm{Z}(A)^{\times}$
s.t.\
\[ \psi^{-1} = a \otimes 1 = 1 \otimes b^{\circ} \in
\opn{End}_{\cat{D}(\cat{Mod} A^{\mrm{e}})}(R_{2}) . \]
So
$\phi = (a \otimes 1) \phi' = (1 \otimes b^{\circ}) \phi'$.
Because $\rho_{2}$ and $\psi$ are $A^{\mrm{e}}$-linear we get
\[ \begin{aligned}
(\phi \otimes \phi) \rho_{1} & =
(a \otimes b^{\circ}) (\phi' \otimes \phi') \rho_{1} = \\
& = (a \otimes b^{\circ}) \rho_{2} \psi \phi' \\
& = \rho_{2} (a \otimes b^{\circ}) \psi \phi' \\
& = \rho_{2} ((a \otimes 1) \psi) (1 \otimes b^{\circ}) \phi' \\
& = \rho_{2} \phi .
\end{aligned} \]
In other words the diagram (\ref{eqn1}) is commutative.
If $\phi''$ also makes (\ref{eqn1}) commutative, then
writing
$\phi'' = (c \otimes 1) \phi$
with $c \in \mrm{Z}(A)^{\times}$, the same computation
shows that $c = 1$.

It remains to produce $\phi'$. Consider the complexes
$T := \mrm{R} \opn{Hom}_{A}(R_{1}, R_{2})$
and
$T^{\circ} := \mrm{R} \opn{Hom}_{A^{\circ}}(R_{1},
R_{2})$.
Then by Theorem \ref{thm4.1}
\[ R_{2}
\cong T^{\circ} \otimes^{\mrm{L}}_{A} R_{1}
\cong R_{1} \otimes^{\mrm{L}}_{A} T . \]
Now using $\rho_{1}$  and $\rho_{2}$ we obtain isomorphisms
in $\cat{D}(\cat{Mod}(A^{\mrm{e}}))$
\[ \begin{aligned}
R_{2}
& \cong \mrm{R} \opn{Hom}_{A^{\mrm{e}}}
(A, R_{2} \otimes R_{2}) \\
& \cong \mrm{R} \opn{Hom}_{A^{\mrm{e}}}
(A, (R_{1} \otimes^{\mrm{L}}_{A} T)
\otimes (T^{\circ} \otimes^{\mrm{L}}_{A} R_{1})) \\
& \cong \mrm{R} \opn{Hom}_{A^{\mrm{e}}}
(A, R_{1} \otimes R_{1})
\otimes^{\mrm{L}}_{A^{\mrm{e}}} (T \otimes T^{\circ}) \\
& \cong R_{1} \otimes^{\mrm{L}}_{A^{\mrm{e}}}
(T \otimes T^{\circ}) \\
& \cong T^{\circ} \otimes^{\mrm{L}}_{A}
R_{1} \otimes^{\mrm{L}}_{A} T \\
& \cong R_{2} \otimes^{\mrm{L}}_{A} T
\end{aligned} \]
so again by Theorem \ref{thm4.1}, $T \cong A$.
\end{proof}

Usually we will leave the isomorphism $\rho$ implicit, and just speak
of a rigid dualizing complex $R$.

\begin{lem} \label{lem5.1}
Suppose $A$ is commutative, integral of dimension $n$ and smooth over
$k$. Then $\Omega^{n}_{A / k}[n]$ is a rigid dualizing complex.
\end{lem}

\begin{proof}
There is a natural isomorphism
$\Omega^{n}_{A / k} \otimes \Omega^{n}_{A / k} \cong
\Omega^{2n}_{A^{\mrm{e}} / k}$ by wedge product.
By \cite{RD} Prop.\ III.8.4 we get a natural isomorphism
\[ \rho : \Omega^{n}_{A / k}[n] \iso
\mrm{R} \opn{Hom}_{A^{\mrm{e}}}(A, \Omega^{2n}_{A^{\mrm{e}} / k}[2n])
. \]
\end{proof}

\begin{rem}
Observe that this $\rho$ is actually the fundamental class of the
diagonal
$X \inj X \times X$, $X = \opn{Spec} A$.
Locally there are generators $a_{1}, \ldots, a_{n}$ for
$\opn{Ker}(A^{\mrm{e}} \ar A)$, and then $\rho$ is given by
the generalized fraction
$\gfrac{\mrm{d} a_{1} \wedge \cdots \wedge \mrm{d} a_{n}}{
a_{1} \cdots a_{n}}$.
\end{rem}

\begin{rem}
J.\ Lipman (in unpublished notes) studied the canonical isomorphism
\begin{equation} \label{eqn5.2}
f^{!} \mcal{O}_{Y} \cong \mrm{R} \mcal{H}om_{X \times_{Y} X}
(\mcal{O}_{X}, f^{!} \mcal{O}_{Y} \boxtimes f^{!} \mcal{O}_{Y})
\end{equation}
where $f : X \ar Y$ is a flat morphism of schemes, in connection
with the relative fundamental class of $f$. When $Y$ is a Gorenstein
scheme, $\mcal{R} := f^{!} \mcal{O}_{Y}$ is a dualizing complex
on $X$. This generalizes Lemma \ref{lem5.1}.
\end{rem}

A ring homomorphism $A \ar B$ is called {\em finite} if $B$
is a finitely generated left and right $A$-module.

\begin{prop} \label{prop5.1}
Suppose $A$ is finite over its center and finitely generated as
$k$-algebra.
Then $A$ has a rigid dualizing complex.
\end{prop}

\begin{proof}
Choose a finite centralizing homomorphism $C \ar A$, with
$C = k[t_{1}, \ldots, t_{n}]$ a commutative polynomial algebra. Let
$R_{C} := \Omega^{n}_{C / k}[n]$, with $\rho_{C}$ as in Lemma
\ref{lem5.1}. Define
$R_{A} := \mrm{R} \opn{Hom}_{C}(A, R_{C})$,
which by \cite{Ye} Proposition 5.2 is a dualizing complex over 
$A$.
One has 
\[ \begin{aligned}
R_{A} \otimes R_{A} & =
\mrm{R} \opn{Hom}_{C}(A, R_{C}) \otimes
\mrm{R} \opn{Hom}_{C}(A, R_{C}) \\
& = \mrm{R} \opn{Hom}_{C^{\mrm{e}}}(A^{\mrm{e}},
R_{C} \otimes R_{C}) .
\end{aligned} \]
Next using $\rho_{C}$ we obtain an isomorphism
\[ \begin{aligned}
\mrm{R} \opn{Hom}_{A^{\mrm{e}}}(A,
R_{A} \otimes R_{A}) & \cong
\mrm{R} \opn{Hom}_{A^{\mrm{e}}}(A,
\mrm{R} \opn{Hom}_{C^{\mrm{e}}}(A^{\mrm{e}},
R_{C} \otimes R_{C})) \\
& \cong \mrm{R} \opn{Hom}_{C^{\mrm{e}}}(A,
R_{C} \otimes R_{C}) \\
& \cong \mrm{R} \opn{Hom}_{C}(A,
\mrm{R} \opn{Hom}_{C^{\mrm{e}}}(C,
R_{C} \otimes R_{C})) \\
& \cong \mrm{R} \opn{Hom}_{C}(A, R_{C}) \\
& \cong R_{A}
\end{aligned} \]
which we label $\rho_{A}$.
\end{proof}

\begin{prop} \label{prop5.2}
Let $A \ar B$ be a finite homomorphism of $k$-algebras, and suppose
$(R_{A}, \rho_{A})$ and $(R_{B}, \rho_{B})$ are
rigid dualizing complexes.
Assume that for some commutative finitely generated $k$-algebra $C$
there exists a homomorphism $C \ar A$, which makes $A$ and $B$
finite $C$-algebras. Then there is a canonical morphism
$\opn{Tr}_{B / A} : R_{B} \ar R_{A}$
in $\cat{D}(\cat{Mod}(A^{\mrm{e}}))$.
\end{prop}

\begin{proof}
Choose such a homomorphism $C \ar A$, and
pick a rigid dualizing complex $(R_{C}, \rho_{C})$.
By Proposition \ref{prop5.1} Theorem \ref{thm5.1} there are unique
isomorphisms
$R_{A} \cong \mrm{R} \opn{Hom}_{C}(A, R_{C})$
and
$R_{B} \cong \mrm{R} \opn{Hom}_{C}(B, R_{C})$.
We obtain $\opn{Tr}_{B / A}$ by applying
$\mrm{R} \opn{Hom}_{C}(-, R_{C})$
to the morphism
$A \ar B$ in $\cat{D}(\cat{Mod} A^{\mrm{e}})$.
This is independent of $C$ by Theorem \ref{thm5.1}.
\end{proof}

\begin{rem}
These results are interesting even for $A$ commutative.
For instance, if $A,B$ are integral of dimension $n$ and smooth
over $k$, and if $A \ar B$ is a finite homomorphism,
then we obtain
$\opn{Tr}_{B / A} : \Omega^{n}_{B / k}[n] \ar \Omega^{n}_{A / k}[n]$.
This trace coincides with the trace of \cite{RD}.
If $A \ar B$ is also \'{e}tale
then
$\Omega^{n}_{B / k} \cong B \otimes_{A} \Omega^{n}_{A / k}$,
and $\opn{Tr}_{B / A}$ is induced from
$B \ar \opn{End}_{A}(B) \xrightarrow{\opn{Tr}} A$.
\end{rem}

Derived equivalent algebras have ``the same'' dualizing complexes:

\begin{prop}
Let $A$ and $B$ be noetherian $k$-algebras,
$R \in \cat{D}(\cat{Mod} A^{\mrm{e}})$ a dualizing complex,
and
$T \in \cat{D}(\cat{Mod} (B \otimes A^{\circ}))$
a tilting complex.
Then
\[ R^{T} := T \otimes^{\mrm{L}}_{A} R \otimes^{\mrm{L}}_{A} T^{\vee}
\in \cat{D}(\cat{Mod} B^{\mrm{e}}) \]
is a dualizing complex.
If in addition $(R, \rho)$ is a rigid dualizing complex, then
$(R^{T}, \rho^{T})$ is rigid, where $\rho^{T}$ is induced naturally
by $\rho$.
\end{prop}

\begin{proof}
Since for any
$M \in \cat{D}^{\mrm{b}}(\cat{Mod} B)$
we have
\[ \mrm{R} \opn{Hom}_{B}(M, R^{T}) \cong
\mrm{R} \opn{Hom}_{A}(T^{\vee} \otimes^{\mrm{L}}_{B} M, R)
\otimes^{\mrm{L}}_{A} T^{\vee} \]
etc.\ it follows that $R^{T}$ is dualizing.

In the rigid situation, first note that
$R^{T} \cong (T \otimes T^{\vee}) \otimes^{\mrm{L}}_{A^{\mrm{e}}} R$,
and
$T \otimes T^{\vee} \in
\cat{D}(\cat{Mod} (B^{\mrm{e}} \otimes (A^{\mrm{e}})^{\circ}))$
is a tilting complex. Using the isomorphism $\rho$ we obtain
\[ \begin{aligned}
\mrm{R} \opn{Hom}_{B^{\mrm{e}}}(B,
R^{T} \otimes R^{T}) & \cong
\mrm{R} \opn{Hom}_{B^{\mrm{e}}}(B,
(T \otimes^{\mrm{L}}_{A} R) \otimes
(R \otimes^{\mrm{L}}_{A} T^{\vee})) \otimes^{\mrm{L}}_{A^{\mrm{e}}}
(T^{\vee} \otimes T) \\
& \cong \mrm{R} \opn{Hom}_{B^{\mrm{e}}}(B,
(T \otimes T^{\vee}) \otimes^{\mrm{L}}_{A^{\mrm{e}}} (R \otimes R))
\otimes^{\mrm{L}}_{A^{\mrm{e}}} (T^{\vee} \otimes T) \\
& \cong \mrm{R} \opn{Hom}_{A^{\mrm{e}}}(A,
R \otimes R)
\otimes^{\mrm{L}}_{A^{\mrm{e}}} (T^{\vee} \otimes T) \\
& \cong R \otimes^{\mrm{L}}_{A^{\mrm{e}}} (T^{\vee} \otimes T) \\
& \cong R^{T} .
\end{aligned} \]
This determines $\rho^{T}$.
\end{proof}

The next proposition generalizes \cite{VdB} Proposition 8.4, which
gives a formula for the rigid dualizing complex $R$ when $A$ is a
Gorenstein algebra and $R \cong L[n]$ for an invertible bimodule $L$.

\begin{prop} \label{prop5.5}
Suppose $A$ is a Gorenstein algebra and $R$ is a rigid dualizing
complex. Then $R$ is a tilting complex and
\[ R^{\vee} = \mrm{R} \opn{Hom}_{A}(R, A)
\cong \mrm{R} \opn{Hom}_{A^{\mrm{e}}}(A, A^{\mrm{e}}) \in
\cat{D}(\cat{Mod} A^{\mrm{e}}) . \]
\end{prop}

\begin{proof}
$R$ is tilting by Theorem \ref{thm4.1}.
Then it is a straightforward calculation:
\[ \begin{aligned}
R & \cong \mrm{R} \opn{Hom}_{A^{\mrm{e}}}(A, R \otimes R) \\
& \cong \mrm{R} \opn{Hom}_{A^{\mrm{e}}}(A, A^{\mrm{e}})
\otimes_{A^{\mrm{e}}}^{\mrm{L}} (R \otimes R) \\
& \cong R \otimes_{A}^{\mrm{L}}
\mrm{R} \opn{Hom}_{A^{\mrm{e}}}(A, A^{\mrm{e}})
\otimes_{A}^{\mrm{L}} R
\end{aligned} \]
so applying
$R^{\vee} \otimes_{A^{\mrm{e}}}^{\mrm{L}} -$
and then
$- \otimes_{A^{\mrm{e}}}^{\mrm{L}} R^{\vee}$
we get what we want.
\end{proof}

\section{Finite $k$-Algebras}

In this section $A$ is a finite $k$-algebra.
We write $M^{*} = D M := \opn{Hom}_{k}(M, k)$ for an $A$-module $M$.
The bimodule $A^{*}$ is then injective on both sides, and
$M^{*} \cong \opn{Hom}_{A}(M, A^{*})$
for any $M \in \msf{D}(\cat{Mod} A)$.

\begin{prop} \label{prop6.1}
\begin{enumerate}
\item $A^{*}$ is a rigid dualizing complex over $A$.
\item $T \in \msf{D}^{\mrm{b}}(\cat{Mod} A^{\mrm{e}})$
is a tilting complex if{f} $T^{*}$ is a dualizing complex.
\item $A$ is a Gorenstein algebra if{f} $A^{*}$ is a tilting 
complex. In this case,
\[ A^{*} \otimes^{\mrm{L}}_{A} M \cong \mrm{R} 
\opn{Hom}_{A}(M, A)^{*} \]
for any
$M \in \msf{D}^{-}_{\mrm{f}}(\cat{Mod} A)$.
\end{enumerate}
\end{prop}

\begin{proof}
1.\ By the proof of Proposition \ref{prop5.1}.
\newline
2.\ Use the duality $D$ (cf.\ Proposition \ref{prop4.1}).
\newline
3.\ Since $A$ is a Gorenstein algebra if{f} $R = A$ is a dualizing
complex, this is a consequence of part 2.
Using a projective resolution of $M$ we get a functorial
morphism
\[ A^{*} \otimes^{\mrm{L}}_{A} M \ar
\opn{Hom}_{A^{\circ}}(\mrm{R} \opn{Hom}_{A}(M, A), A^{*}) . \]
By way-out arguments it suffices to check that this is an isomorphism
for $M = A$, which is clear.
\end{proof}

\begin{rem}
When the dualizing complex $R$ is a single bimodule in degree $0$,
it is called a {\em cotilting module} in the literature.
The name is justified by part 2 of the proposition
(and cf.\ Theorem \ref{thm4.1}).
\end{rem}

\begin{rem}
The derived functor $A^{*} \otimes^{\mrm{L}}_{A} -$ is discussed
in \cite{Ha} and in \cite{Ri2} Section 5.
If $A$ is a hereditary algebra then by Proposition \ref{prop6.1}(3)
we have
\[ \mrm{H}^{-1} (A^{*} \otimes^{\mrm{L}}_{A} M) \cong
\opn{Ext}^{1}_{A}(M, A)^{*} \cong D \opn{Tr} M  \]
for every $M \in \cat{Mod}_{\mrm{f}} A$. Here $D \opn{Tr}$
is the `dual of the transpose' functor of \cite{ARS} Chapter IV,
which induces the translation function in the
{\em Auslander-Reiten quiver} of $A$.
\end{rem}

Now assume $A$ has finite global dimension.
Let $S_{1}, \ldots, S_{n}$ be a complete set of nonisomorphic
simple $A$-modules, and let $P_{1}, \ldots, P_{n}$ (resp.\
$I_{1}, \ldots, I_{n}$) be the corresponding indecomposable
projective (resp.\ injective) modules. Then the Grothen\-dieck group
$\opn{K}_{0}(A) = \opn{K}_{0}(\cat{Mod}_{\mrm{f}} A)$ is a free
$\mbb{Z}$-module with basis either of the sets
$\{ [S_{i}] \}_{i = 1}^{n}$, $\{ [P_{i}] \}$ or $\{ [I_{i}] \}$.
The {\em Coxeter transformation} $c \in \opn{Aut}(\mrm{K}_{0}(A))$
is defined by
$c([P_{i}]) := -[I_{i}]$
(see \cite{ARS} Section VIII.2).

In Proposition \ref{prop3.5} we defined the representation
$\chi_{0} : \opn{DPic}(A) \ar \opn{Aut}(\mrm{K}_{0}(A))$.
Denote by $t$ the class of $A^{*}$ in $\opn{DPic}(A)$.

\begin{prop} \label{prop6.2}
$\chi_{0}(t) = -c$.
\end{prop}

\begin{proof}
This follows from Proposition \ref{prop6.1}(3) and \cite{ARS}
Proposition VIII.2.2 (a).
\end{proof}

For the remainder of the section we shall examine the algebra
\[ A = \left[ \begin{matrix}
k & k \\
0 & k
\end{matrix} \right] . \]
(This was suggested by T.\ Stafford.)
Observe that $A$ is the smallest $k$-algebra which is neither
commutative nor local, so Propositions \ref{prop3.3} and
\ref{prop3.4} do not apply. In the classification by Dynkin quivers
(diagrams), the algebra $A$ corresponds to the quiver
$\Delta = A_{2} = (\bullet \longrightarrow \bullet)$.
That is, $A \cong k \Delta$, the path algebra of $\Delta$.

Let $P_{1}, P_{2}$ (resp.\ $S_{1}, S_{2}$) be the projective
(resp.\ simple) $A$-modules
\[ P_{1} = S_{1} :=
\left[ \begin{matrix} k \\ 0 \end{matrix} \right]\ ;\
P_{2} := \left[ \begin{matrix} k \\ k \end{matrix} \right]\ ;\
S_{2} := \left[ \begin{matrix} 0 \\ k \end{matrix} \right] , \]
so that
$A = P_{1} \oplus P_{2}$ as $A$-modules.

\begin{prop} \label{prop6.3}
\begin{enumerate}
\item $\opn{Pic}(A) = 1$.
\item There are isomorphisms in $\msf{D}(\cat{Mod} A)$:
\[ \begin{aligned}
A^{*} \otimes^{\mrm{L}}_{A} S_{1} & \cong P_{2} \\[1mm]
A^{*} \otimes^{\mrm{L}}_{A} P_{2} & \cong S_{2} \\[1mm]
A^{*} \otimes^{\mrm{L}}_{A} S_{2} & \cong S_{1}[1] .
\end{aligned} \]
\item There is an isomorphism in $\msf{D}(\cat{Mod} A^{\mrm{e}})$:
\[ A^{*} \otimes^{\mrm{L}}_{A} A^{*} \otimes^{\mrm{L}}_{A} A^{*}
\cong A[1] . \]
\end{enumerate}
\end{prop}

\begin{proof}
1.\ First note that the indecomposable projective modules
$P_{1}$ and $P_{2}$ have different lengths.
So if $L$ is an invertible bimodule we must have
$L \otimes_{A} P_{1} \cong P_{1}$ and
$L \otimes_{A} P_{2} \cong P_{2}$.
Therefore
$L \cong A$ as $A$-modules.
According to Lemma \ref{lem3.1}(2) we get
$L \cong A_{\sigma}$ as bimodules, for some $\sigma \in \opn{Aut}(A)$.
But one sees that any such $\sigma$ is conjugation by a matrix
$\left[ \begin{smallmatrix}
a & b \\ 0 & 1 \end{smallmatrix} \right]$,
so
$A_{\sigma} \cong A$ as bimodules and
$\opn{Pic}(A) \cong \opn{Out}(A) = 1$.
\newline
2.\ A straightforward calculation using the isomorphism
of $A^{\mrm{e}}$-modules
\[ A^{*} \cong
\left[ \begin{smallmatrix}
k & 0 \\
k & k
\end{smallmatrix} \right] =
\left[ \begin{smallmatrix}
k & k \\
k & k
\end{smallmatrix} \right] /
\left[ \begin{smallmatrix}
0 & k \\
0 & 0
\end{smallmatrix} \right] \]
induced by the trace pairing on
$\mrm{M}_{2}(k) = \left[ \begin{smallmatrix}
k & k \\
k & k
\end{smallmatrix} \right]$.
\newline
3.\ By part 2 we obtain this isomorphism in
$\msf{D}(\cat{Mod} A)$. Now apply Proposition \ref{prop2.1},
Lemma \ref{lem3.1} and part 1 above.
\end{proof}

As before denote by $s$ the class of $A[1]$ in $\opn{DPic}(A)$.
The action of $s$ on $\msf{D}(\cat{Mod} A)$ is by
a shift in degree, and the subgroup $\langle s \rangle$
is then isomorphic to $\mbb{Z}$.
Part 3 of the proposition gives the remarkable fact:

\begin{cor}
$t^{3} = s$.
\end{cor}

In terms of the representation $\chi_{0}$ and the basis
$\{ [S_{1}], [S_{2}] \}$ of $\mrm{K}_{0}(A)$ we get
\[
\chi_{0}(s) =
\left[ \begin{matrix}
-1 & 0 \\
0  & -1
\end{matrix} \right]
\qquad
\chi_{0}(t) =
\left[ \begin{matrix}
1   & 1 \\
-1  & 0
\end{matrix} \right] .
\]

\begin{rem}
These results were extended by E.\ Kreines to upper triangular
$n \times n$ matrix rings, $n \geq 2$ (see the Appendix).
In particular she showed that
$t^{n+1} = s^{n-1}$.
This is in agreement with the fact that the order of the Coxeter
transformation $c$ is $n + 1$, cf.\ \cite{ARS} p.\ 289.
\end{rem}

\begin{prob} \label{prob6.1}
Let $A$ be an indecomposable, elementary, hereditary $k$-algebra of
finite representation type.
What is the structure of the group $\opn{DPic}(A)$?
Is it true that
$\opn{DPic}(A) \cong \mbb{Z}$ with generator $t$?
What is the structure of the rings $\opn{K}^{0}(A)$
and $\opn{DK}^{0}(A)$?
How do $\opn{DPic}(A)$ and $\opn{DK}^{0}(A)$ fit in with other
invariants of $A$?
\end{prob}

\appendix
\section{The Algebra of $n \times n$ Upper Triangular Matrices 
\protect\newline
by Elena Kreines}

Let us consider the upper triangular $n \times n$ matrix algebra 
$A$ over a field $k$, where $n \geq 2$. 
Let $A^{*} := \opn{Hom}_{k}(A,k)$, which is known to be a tilting 
complex, and define the functor
$F: \cat{D}(\cat{Mod} A) \ar \cat{D}(\cat{Mod} A)$,
$F M := A^{*} \otimes^{\mrm{L}}_{A} M$.

\begin{thm}
There is an isomorphism
\[ F^{n+1} A =
\underbrace{A^{*} \otimes^{\mrm{L}}_{A}
\ldots \otimes^{\mrm{L}}_{A} A^{*}}_{n+1} 
\otimes^{\mrm{L}}_{A} A \cong A[n-1] \]
in $\cat{D}(\cat{Mod} A)$.
\end{thm}

The proof of the theorem appears at the end of the appendix. 

\begin{cor} 
We get an isomorphism
\[ \underbrace{A^{*} \otimes^{\mrm{L}}_{A} \ldots 
\otimes^{\mrm{L}}_{A} A^{*}}_{n+1} \cong A[n-1] \]
in $\cat{D}(\cat{Mod} A^{\mrm{e}})$. Hence 
$t^{n + 1} = s^{n - 1}$ in $\opn{DPic}(A)$.
\end{cor}

\begin{proof}[Proof of the corollary]
By \cite{BK}, 
$\opn{Aut}(A)= \opn{Inn}(A)$, and thus we can use the proof  of
Proposition \ref{prop6.3}(3). 
\end{proof}

Let $\opn{M}_{n}(k)$ denote the full matrix algebra, and let
$\mfrak{r} \subset A$ be the ideal of strictly upper triangular 
matrices. Then the trace pairing on $\opn{M}_{n}(k)$ identifies
$A^{*} \cong \opn{M}_{n}(k) / \mfrak{r}$
as $A$-bimodules.

For $1 \leq i \leq j \leq n$ let $I^{i}_{j}$ be the $A$-module 
represented as a column
\[ I^{i}_{j} :=
\left[ \begin{array}{c}
0 \\ : \\ 0 \\ k \\ : \\ k \\ 0 \\ : \\ 0
\end{array} \right]
\begin{array}{c}
\blnk{2mm} \\ \blnk{2mm} \\ \blnk{2mm} \\ \leftarrow i \\ 
\blnk{2mm} \\ \leftarrow j \\ \blnk{2mm} \\ \blnk{2mm} \\ 
\blnk{2mm} \end{array} \]
The left action of $A$ on $I^{i}_{j}$ is as follows. For $i = 1$ 
this is the usual matrix multiplication, and for $i > 1$ we have
$I^{i}_{j} \cong I^{1}_{j} / I^{1}_{i - 1}$. 

We see that $P_{j} := I^{1}_{j}$ is a projective module, and 
$A = \bigoplus^{n}_{j = 1} P_{j}$. 
Also $I_{i}:=I^{i}_{n}$ is an injective module, and 
$A^{*} = \bigoplus^{n}_{i = 1} I_{i}$.
The module $S_{i} :=I^{i}_{i}$ is simple.

For the proof of the theorem we need two lemmas.

\begin{lem} \label{lem.ap.1} 
For $i = 1, \ldots, n$ we have
$ FP_{i} \cong I^{i}_{n}$.
\end{lem}

\begin{proof}
Since the module $P_i$ is projective we have 
$A^{*} \otimes^{\mrm{L}}_{A} P_{i} = A^{*} \otimes_{A} P_{i}$. 
By tensoring the short exact sequence  
\[ 0 \ar \bigoplus_{\substack{j = 1, \ldots , n \\ j \neq i}} 
P_{j} \ar A \ar  P_{i} \ar 0 \]
with the module $A^{*}$, and noting that 
$A^{*} \cdot P_{j} = I^{j}_{n} \subset A^{*}$
(where we view $P_{j} \subset A$ as a left ideal), we obtain
\[ F P_{i} = A^{*} \otimes_{A} P_{i} \cong 
\frac{A^{*}}{A^{*} \cdot (\bigoplus_{j \neq i} P_{j})}
= \frac{A^{*}}{\bigoplus_{j \neq i} I^{j}_{n}}
\cong I^{i}_{n} . \]
\end{proof}

\begin{lem} 
If $i > 1$ then $F I^{i}_{j} \cong I^{i-1}_{j-1}[1]$.
\end{lem}

\begin{proof}
The module $I^i_j$ is not projective.
A projective resolution for this module is
the short exact sequence 
\[ 0 \ar P_{i - 1} \ar P_{j} \ar I^{i}_{j} \ar 0 . \]
By tensoring this sequence with the module $A^{*}$ and using Lemma
\ref{lem.ap.1} we obtain the exact sequence
\[  I^{i - 1}_{n} \xrightarrow{\phi} I^{j}_{n} \ar 
A^{*} \otimes_{A} I^{i}_{j} \ar 0 . \]

Let us denote by $M_{1} \subset A$ the set of matrices whose only 
nonzero entries are in the first row. It is easy to see that
$M_{1} \cdot I^{i}_{j} = 0$ (since $i > 1$) and that 
$A^{*} \cdot M_{1} = A^{*}$. This implies that
$\opn{Coker}(\phi) = A^{*} \otimes_{A} I^{i}_{j} = 0$.
Since $\opn{Ker}(\phi)$
is a submodule of $I^{i - 1}_{n}$ of length $j - i + 1$ 
we must have
$\opn{Ker}(\phi) = I^{i - 1}_{j - 1}$.
\end{proof}

\begin{proof}[Proof of the theorem]
By the two lemmas
\[ \begin{array}{ccccc}
F P_{i} & \cong & I^{i}_{n}  \\[2mm]
F^{2} P_{i} & \cong & I^{i-1}_{n-1}[1] \\
& \vdots \\
F^{i} P_{i} & \cong & I^{1}_{n + 1 - i}[i - 1] & 
= & P_{n + 1 - i}[i - 1] \\
& \vdots \\ 
F^{n+1} P_{i} & \cong & P_{i}[n-1] .
\end{array} \]
But $A = \bigoplus_{i=1}^{n} P_{i}$ as $A$-modules, and hence
$F^{n+1} A \cong A[n-1]$ as claimed.
\end{proof}


\end{document}